\documentclass[proc]{edpsmath}
\usepackage{graphicx}
\usepackage{color}
\usepackage{amsmath, amssymb}
\usepackage[english]{babel}
\usepackage{placeins}
\usepackage[utf8]{inputenc}
\usepackage{hyperref}
\usepackage{caption,subfigure}

\usepackage{color}				% Pour la couleur
\usepackage[colorinlistoftodos,prependcaption]{todonotes}	% Pour les notes (pour désactiver: \usepackage[disable]{todonotes})

\definecolor{LGcolor}{RGB}{190,23,31}									

\definecolor{AVcolor}{RGB}{241,195,49}									

\definecolor{ALcolor}{RGB}{49,184,114}									

\definecolor{DMcolor}{RGB}{44,141,191}

\usepackage[bbgreekl]{mathbbol}  %  Pour afficher les tenseurs avec deux barres verticales
\usepackage{lmodern}
\usepackage{anyfontsize}

\makeatletter
\def\amsbb{\use@mathgroup \M@U \symAMSb}
\makeatother
\newcommand{\beq}[1]{\begin{equation}\label{#1}}
\newcommand{\eeq}{\end{equation}}
\newcommand{\bitem}{\begin{itemize}}
\newcommand{\eitem}{\end{itemize}}
\newcommand{\benum}{\begin{enumerate}}
\newcommand{\eenum}{\end{enumerate}}

\newcommand{\G}{\boldsymbol{G}}

\newcommand{\dzeta}{\boldsymbol{\zeta}}

\newcommand{\bx}{\pmb{x}}
\newcommand{\bz}{\pmb{z}}
\newcommand{\bZ}{\pmb{Z}}

\newcommand{\bX}{\pmb{X}}
\newcommand{\bc}{\pmb{c}}
\newcommand{\bC}{\pmb{C}}

\newcommand{\bW}{\pmb{W}}
\newcommand{\nul}{\boldsymbol{0}}
\newcommand{\bH}{\boldsymbol{H}}

\newcommand{\bF}{\pmb{F}}
\newcommand{\Fext}{\bF_{\text{ext}}}

\newcommand{\bn}{\pmb{n}}

\newcommand{\bu}{\pmb{u}}
\newcommand{\uf}{\overline{\bu_f}}

\newcommand{\bq}{\pmb{q}}

\newcommand{\R}{\mathbb{R}}

\renewcommand{\P}{\mathbb{P}}
\newcommand{\ds}{\displaystyle}
\newcommand{\dd}[2]{\frac{\partial #1}{\partial #2}}

\newcommand{\dt}[1]{\dfrac{d #1}{dt}}

\renewcommand{\div}{\boldsymbol{\nabla}\cdot}
\newcommand{\Div}[1]{\div\left( #1 \right)}
\newcommand{\grad}{\boldsymbol{\nabla}}

\newcommand{\eps}{\varepsilon}
\newcommand{\F}{\boldsymbol{F}}

\newcommand{\tw}{\textwidth}

\newcommand{\I}{\amsbb{I}}

\newcommand{\aire}[1]{\left|#1\right|}
\newcommand{\moy}[1]{\overline{#1}} %left<#1\right>}

\renewcommand{\H}{\mathcal{H}}
\newcommand{\feq}{f_{eq}}

\newcommand{\Kn}{\text{Kn}}

\begin{document}

%%-----------------------------
%%      the top matter
%%-----------------------------
\title{Statistical and probabilistic modeling of a cloud of particles coupled with a turbulent fluid}
\author{Ludovic Gouden\`ege}\address{Fédération de Mathématiques de CentraleSupélec, CNRS FR-3487, CentraleSup\'elec, Universit\'e Paris-Saclay, 9 rue Joliot Curie, 91190 Gif-sur-Yvette cedex, France}
\author{Adam Larat}
\sameaddress{1, 6}
\secondaddress{Univ. Grenoble Alpes, CNRS, Grenoble INP, LJK, 38000 Grenoble, France}
%Laboratoire EM2C, CNRS, CentraleSup\'elec, Universit\'e Paris-Saclay, 3 rue Joliot Curie, 91192 Gif-sur-Yvette cedex, France.}
\author{Julie Llobell}\address{{Universit\'e C\^ote d'Azur, Inria,  CNRS, LJAD, Parc~Valrose, F-06108 
Nice, France}}
\author{Marc Massot}\address{CMAP, \'Ecole Polytechnique, Route de Saclay, 91128 Palaiseau cedex, France}
\author{David~Mercier}\sameaddress{6}
\author{Olivier Thomine}\address{CEA/DEN/DANS, 91191 Gif-sur-Yvette cedex, France}
\author{Aymeric Vi\'e}\sameaddress{1}\secondaddress{Laboratoire EM2C, CNRS, CentraleSup\'elec, Universit\'e Paris-Saclay, 3 rue Joliot Curie, 91192 Gif-sur-Yvette cedex, France}

\begin{abstract} 
  This paper exposes a novel exploratory formalism, which end goal is 
  the numerical simulation of the dynamics of 
  a cloud of particles weakly or strongly coupled with a turbulent fluid. 
  Given the large panel of expertise of the list of authors, 
  the content of this paper scans a wide range of connex notions, 
  from the physics of turbulence to the rigorous definition of stochastic processes. 
  Our approach is %, starting from the formalism of Pope \cite{pope2010self}, which sets an interesting rigorous mathematical framework for the physical description of turbulent fluids, 
  to develop reduced-order models for 
  the dynamics of both carrying and carried phases which remain consistant 
  within this formalism, and to set up a numerical process to validate these models. 
  The novelties of this paper lie in the gathering of a large panel of mathematical 
  and physical definitions and results within a common framework and an 
  agreed vocabulary (sections \ref{sec:MicroMeso} and \ref{sec:Macro}), 
  and in some preliminary results and achievements within this context, section \ref{sec:Numerics}. 
  While the first three sections have been simplified to the context of a gas field providing 
  that the disperse phase only retrieves energy through drag, 
  the fourth section opens this study to the more complex situation when the disperse phase 
  interacts with the continuous phase as well, in an energy conservative manner. This will 
  allow us to expose the perspectives of the project and to conclude.  
\end{abstract}
%
%\begin{resume} \end{resume}
%
%
\maketitle

\tableofcontents

%%-----------------------------
%%      your text
%%-----------------------------
\section*{Introduction}

Many applications involve the transport of a disperse phase (particles, droplets, bubbles) coupled with a fluid: spray combustion, fluidized beds, soot dynamics... % \cite{Zamansky16,Vie13}. 
In the standard case, the evolution of the carrier phase can be
described by a deterministic system of equations such as the Navier-Stokes equations.
However, in the strongly coupled case the evolution equations are unclosed due to the
exchange term with the particles \cite{Fox14,Emre14}. 
Often, models proposed in the literature only consider the influence of the carrier fluid
on the disperse phase and neglect its retroactive consequences, or, at best, limit it to a global balance
between the two phases \cite{minier2016statistical}. In particular, these simplifying hypotheses 
allow to decouple the inaccuracies coming from the approximate resolution of each phase.

%%%Adam%%% In order to compute the evolution of particles in turbulence, different levels of simulation are possible.
%%%Adam%%% The more detailled one is achieved by fully resolving the evolution of the particle and the dynamic of the flow around it. 
%%%Adam%%% When the particle is spherical and its size is of the order of Kolmogorov length scale, it is possible to consider that the particle is a point in the fluid with no physical extension. 
%%%Adam%%% Furthermore, when its density is bigger than that of the fluid, it is also possible to only consider the drag force. 
%%%Adam%%% When the particle relative velocity is small, the drag force can be linearized. It is then called Stokes drag law. 
%%%Adam%%% In the following, for the numerical simulations, we will only consider Stokes drag law to be the most predominant force acting on the particle evolution. 
%%%Adam%%% Thus particle evolution is thus given by equation \eqref{eq:partEvoDNS} . For the moment there is no coupling with the flow $\boldsymbol{v}_f$.

But, one of the main difficulties in the derivation of a consistent model for the strongly coupled 
evolution of a cloud of particles within a turbulent flow, is that inaccuracies arise 
both from the chaotic behavior of the fluid \cite{Pope01,Zaichik09,Weinan01}, and from the initial properties 
of the particles, such as their starting positions and velocities. Therefore, the proper level of ``modeling'' 
consists in making consistant assumptions about the %statistical and probabilistic 
properties of the stochastic processes involved in the global dynamics of both phases.
Even if some advances have been made in the field \cite{Gorokhovski14}, the problematic is far
from being closed. 

In order to better understand the coupling of the inaccuracies coming from both phases, 
we split the construction of the fluid dynamic model into four main steps, corresponding 
to four spatial levels of modeling. Step-by-step, we then express some links between these levels, in order to 
better understand the influence of the small scales on the highest level of modeling. 
Here, one has to understand that this hierarchy of points of view is worth both for the carrier fluid 
and for the disperse phase. Simply, each of the passing to the limit between each levels does not 
occur at the same scale for the two phases. Although the carrier fluid is made of nanometric particles, while the dispersed particles 
seldom reach a micrometrical size, the description of each phase starts at the microscopic level (or molecular level). 
From there, one can reach reduced order large-scale models rather continuously, by first looking at 
an intermediate \textit{mesoscopic} scale, dealing with the law on the presence of the 
microscopic phase (\textit{e.g.} the Boltzmann equation), 
and then consider close-to-equilibrium regimes that we will call the \textit{macroscopic} scale
(\textit{e.g.} Euler or Navier-Stokes equations). These four different levels of modeling are 
sketched level-by-level in the following items list:
\begin{itemize}
  \item \underline{\textbf{Microscopic:}} at the scale of atomes, molecules or particules. 
    Generally speaking, one may say "\textit{at the scale of the indivisible}". 
    The medium is here modeled by a very large number of ODEs. 
  \item \underline{\textbf{Macroscopic:}} at the scale of the continuum. 
    Fluids (liquid, gaz, spray,\dots) are now seen as a continuous medium. 
    It is modeled by a system of PDEs. 
  \item \underline{\textbf{Mesoscopic:}} the transition from the micro to the macro scale 
    necessitates an intermediate scale, called "\textit{mesoscopic}", 
    at which the medium is modeled a statistical manner. 
    At this level, the fluid is modeled by the transport equation of a probability density function (PDF) of particles. 
  \item \underline{\textbf{Reduced-Order:}} despite all the complexity reduction already performed, 
    the simulation of all the macroscopic scales (Direct Numerical Simulation, DNS) 
    is far from being reachable. An additional order reduction is then performed 
    by splitting the solution into a \textit{significant} part
    $\moy{\phi}$ and a \textit{residual} $\phi'$: $\phi = \moy{\phi} + \phi'$. 
    In general, the residual is removed and its action on the 
    resolved part is modeled by a chosen underlying random process.
\end{itemize}

Throughout this paper the term \textit{significant part} is kept general on purpose: it could denote 
one of the numerous choice of decomposition of the macroscopic sought solution into a numerically resolved 
and an unresolved part, see subsection \ref{ssec:ProbaMapPope} for more details. 
To give an insight of historical context, the usual method is traditionally referred to as a Large Eddy Simulation (LES) of the particulate flow, 
which means that only the features of the flow at a scale greater than a characteristic cut-off size are 
computed. The smallest scales, called \textit{subscales}, need to be modeled from the computed variables in both carrying and disperse phases.
As proposed by Pope \cite{pope2010self}, 
we chose to place ourselves in a probabilistic formalism 
where the closure in performed by the definition of a probabilistic process for the residuals.
This closure can be seen as a probabilistic mapping between the reduction of the non-linear terms of the solved macroscopic PDEs and the resolved variables, see subsection \ref{ssec:ProbaMapPope}.
As a consequence, defining a subscale model is equivalent to making a choice for this mapping. This is what we are looking for in this project. 

An ideal model for the numerical simulation of a turbulent flow loaded with dispersed particles would be 
a global reduced-order model for both phases, where the residual part would have to be able to take into account the strong coupling between both phases 
(mass, momentum and energy are exchanged in a bidirectional manner and globally conserved).
We think that the formalism introduced in \cite{pope2010self},
and rapidly sketched in the previous paragraph, is a good starting point.
We also believe that the stochastic model of the unresolved fluctuations has a root at the microscopic level in both phases. 
This is the reason why we then start our exploratory study by considering an idealistic micro/micro modeling 
with additional stochastic processes on both phases, and then try to derive a global large-scale reduced-order model
for the dynamics of the strongly coupled system, which remains reliable, accurate and consistant with the underlying 
micro/micro description of the physical system. 

This paper is divided into four sections. In a first section, we give a statistical description at micro and mesoscale which are the beginning of all macroscopic descriptions, with a theorem in the infinite population limit. It explains the link between a system of a large number of ODEs at the microscopic level and a PDE on a Probability Density Function (PDF) of existence of the particles. 
In section two we describe in a very condensed manner the 
other levels of continuous description, while staying as consistant as possible. 
This leads us to a very general definition of turbulence and to the probabilistic framework for the modeling 
of the subscales in the context described by Pope \cite{pope2010self}.
In particular we explain the derivation of a reduced-order model for the disperse phase only, when the underlying 
carrying continuous gas field is supposed to be perfectly known and is not perturbed by the presence of the 
particles. 
Section three presents a numerical process intended to validate the reduced-order models possibly created within 
this micro/micro to reduced-order context, by looking at the statistics missed by the disperse field when the 
underlying gas velocity field has been reduced (for example filtered). In particular, we show that it seems hard to build a reduced-order turbulent model for the dynamics of a 1D spray, but that the situation improves with higher 
dimensionality. 
Finally, section four opens the discussion on the construction of a consistent model for two-way coupled systems.
This section being preliminary, this will allow us to expose the perspectives of the current project and to conclude the paper.

%\section{Statistical description of turbulent particulate flows}
%\section{Turbulent particulate flows: from microscopic to reduced-order models}
\section{Statistical description of the dynamics of a population: from micro- to meso-scale}
\label{sec:MicroMeso}

This section describes the dynamics of a population at micro and mesoscale. This is the beginning of all work implying complex dynamics of turbulent particules-laden flows.
This gathering represents a real team effort, especially in finding a common vocabulary 
between those of us more physics-oriented and those more used to the theory of probability and of stochastic processes.
As already said in the previous paragraph, what is written here is valid for both carrying and disperse phases, 
only the passing to the limit do not occur at the same scales. 

\subsection{Microscopic scale}

The studied domain $\mathbb{X}\subset \R^{3}$ is filled with a cloud of $N$ identical spherical particles, moving into void or supported by a carrying gas. 
Assuming that the three degrees of freedom in rotation of each particle can be ignored, the dynamics of the system is described by the $6N$ parameters (velocity are in $\R^{3}_{\bC}:=\R^{3}$):
\beq{eq:PhaseLiouville}\bZ(t) = \left(\bX_1(t),\bC_1(t),\dots,\bX_N(t),\bC_N(t)\right)\in \dzeta^N :=\left(\mathbb{X}\times\R^3_{\bC}\right)^N , \eeq 
or equally by the \textit{empirical measure} or \textit{normalized counting measure}: $\mu_t^N[\bZ] = \frac{1}{N} \sum_{i=1}^N \delta_{\bX_i(t)}\delta_{\bC_i(t)}$. 

If the set of particles is immersed within an external field $\G(t,\bX,\bC)$, interacts with itself following a collision kernel $\F(\bX,\bC)$ and each particle is possibly subject to an independent Brownian random process of intensity $\sigma$, 
the phase space \eqref{eq:PhaseLiouville} evolves with the following system of $6N$ ODEs: 
\beq{eq:Liouville}
 \left\{
   \begin{array}{ccl}
     \ds d\bX_i &=& \bC_i(t) dt,\\
     \\
     \ds d\bC_i &=& \G\left(t,\bX_i,\bC_i\right) dt 
                 + \F*\mu_t^N\left(\bX_i,\bC_i\right) dt
                 + \sqrt{2}\sigma d\bW_i(t),
   \end{array}
 \right.
 \quad i=1,\dots,N.
\eeq

Then, given an initial condition $\bZ_0 = \left(\bX_1^0,\bC_1^0,\dots,\bX_N^0,\bC_N^0\right)$, which may  
be deterministic or stochastic, 
%\lgo[inline]{bizarre, déterministe ou pas ? Ici y'a W donc stochastic. }
the empirical measure can be 
indexed by $\bZ_0$: $\mu_t^N[\bZ_0] = \ds\frac{1}{N} \sum_{i=1}^N \delta_{\bX_i(t,\bZ_0)}\delta_{\bC_i(t,\bZ_0)}$, 
so that if $V_{\bX}\times V_{\bC} \subset \dzeta$ is a subset of phase space, 
\[
n_{V_{\bX}\times V_{\bC}} := N.\mu_t^N[\bZ_0] (V_{\bX}\times V_{\bC}) = \sum_{i=1}^{N} \mathbb{1}_{V_{\bX}}(\bX_i(t,\bZ_0))\mathbb{1}_{V_{\bC}}(\bC_i(t,\bZ_0))
\]
is the number of particles from the configuration $\bZ_0$ at time $t=0$, situated within $V_{\bX}$ and with a velocity belonging to $V_{\bC}$ at time $t$. 

\subsection{Mesoscopic scale}

%\bitem
%\item

From now on, the configuration of each particle is denoted by $\bz_i = \left(\bX_i,\bC_i\right)\in \mathbb{X} \times \R^3_{\bC}$, for all $i=1,\dots,N$.
The collision kernel $\F$ simulates the interaction between the particles and it thus seems fair to have 
$\F(-\bz)=\F(\bz)$, which implies $\F(\nul)=\nul$. 

Let us consider that the particles are changeable at initial time, which means that  
their initial distribution $\mu_0^{N} \in \mathbb{R}^{2d}$ is invariant by permutation of the $N$ variables.  
This invariance therefore remains satisfied at any time $t>0$ and in particular, the $N$ particles must follow the same one-particle law in $\mathbb{R}^{2d}$, denoted  $\mu_t^{(1)}$, which is what we are looking for in this subsection.
First, if $A$ is a borelian in $\mathbb{R}^{2d}$,
$$\mathbb{E}[ \mu_t^N(A)]=\mathbb{E} \left[\frac{1}{N} \sum_{i=1}^N \delta_{\bz_i(t)}(A)\right]=\frac{1}{N}\sum_{i=1}^N \mathbb{P}[\bz_i(t)\in A] = \mathbb{P}[\bz_1(t)\in A] = \mu_t^{(1)}(A).$$
Then, we recall $\bZ(t)=\left( \bz_1(t), ...  , \bz_N(t)\right)$, later simply noted $\bZ_t$, 
and we introduce the following function
$$\bH(t,\bZ)= \left( \bc_1, \G(t,\bz_1)+\frac{1}{N}\sum_{j=1}^N \F(\bz_1-\bz_j), \cdots ,  \bc_N, \G(t,\bz_N)+\frac{1}{N}\sum_{j=1}^N \F(\bz_N-\bz_j)\right),$$
and the $2d\times N$ diagonal matrix, denoted $\Sigma$,  with zero ($d$ times) and $\sigma$ ($d$ times), repeated $N$ times along the diagonal. 
Thus, equation \eqref{eq:Liouville} can be rewritten $$d \bZ(t)= \bH(t,\bZ(t))dt+\sqrt{2}\Sigma d\bW(t).$$
For any function $\Phi:(t,\bx) \longmapsto \Phi(t,\bx)$ such that
\begin{equation}\label{hyp_Phi}
t \longmapsto \Phi(t,.) \in \mathcal{C}^1 \text{ and } \bx \longmapsto \Phi(.,\bx) \in \mathcal{C}^{\infty}_c,
\end{equation}
the Itô's formula gives us:
\begin{align*}
\Phi(t,\bZ_t)-\Phi(0,\bZ_0)
&= \sigma^2 \sum_{i=1}^N \int_0^T \Delta_{\bc_i}\Phi(s,\bZ_s)ds + \int_0^T \frac{\partial}{\partial t}\Phi(s,\bZ_s)ds,\\
&+ \sum_{i=1}^N \int_0^T \nabla_{\bx_i}\Phi(s,\bZ_s)\cdot d\bZ_s
 + \sum_{i=1}^N \int_0^T \nabla_{\bc_i}\Phi(s,\bZ_s)\cdot d\bZ_s\\
&= \sigma^2 \sum_{i=1}^N \int_0^T \Delta_{\bc_i}\Phi(s,\bZ_s)ds + \int_0^T \frac{\partial}{\partial t}\Phi(s,\bZ_s)ds \\
&+ \sum_{i=1}^N \int_0^T \nabla_{\bx_i}\Phi(s,\bZ_s) \cdot \bH(s,\bZ_s)ds +\sqrt{2}\sum_{i=1}^N \int_0^T \nabla_{\bx_i}\Phi(s,\bZ_s)\cdot \Sigma d\bW_s\\
&+ \sum_{i=1}^N \int_0^T \nabla_{\bc_i}\Phi(s,\bZ_s) \cdot \bH(s,\bZ_s)ds+\sqrt{2} \sum_{i=1}^N \int_0^T \nabla_{\bc_i}\Phi(s,\bZ_s)\cdot \Sigma d\bW_s,
\end{align*}
where $[\nabla_{\bx_i}\Phi(s,\bZ_{s})\  \cdot\ ]$ and $[\nabla_{\bc_i}\Phi(s,\bZ_{s})\  \cdot\ ] $ 
denote projection operators on the respective lines of $\bx_i$ and $\bc_i$.

Taking the expectancy we get:
\begin{align*}
\mathbb{E}[ \Phi(t,\bZ_t)-\Phi(0,\bZ_0)]
&=\sigma^2 \sum_{i=1}^N \int_0^T \mathbb{E}[ \Delta_{\bc_i}\Phi(s,\bZ_s)]ds + \int_0^T \mathbb{E}\left[ \frac{\partial}{\partial t}\Phi(s,\bZ_s)\right]ds\\
&+\sum_{i=1}^N \int_0^T \mathbb{E}\left[ \nabla_{\bx_i}\Phi(s,\bZ_s) \cdot \bH(s,\bZ_s)\right]ds
+\sum_{i=1}^N \int_0^T \mathbb{E}\left[ \nabla_{\bc_i}\Phi(s,\bZ_s) \cdot \bH(s,\bZ_s)\right]ds.
\end{align*}
Next, we introduce the following linear form on the measures of $\R^{2d}$, defined for any $\Phi$ such as in \eqref{hyp_Phi}:
\[
\langle\mu_t^{(N)},\Phi\rangle=\int_0^T \mathbb{E} [ \Phi(t,\bZ_t) ] dt = \int_0^T \int_{\R^{2d}} \Phi(t,\bz) d\mu_t^{(N)}(\bz) dt.
\]
Here $\mu_t^{(N)}$ is the $N$-joint law followed by the $N$ particules: it is the law followed by $\bZ_t$. 
Using this dual formulation, we can now extend the definition of the partial derivatives to the measures of $\R^{2d}$, and we have:
\begin{align*}
&\int_0^T \mathbb{E}[ \nabla_{\bx_i}\Phi(s,\bZ_s) \cdot \bH(s,\bZ_s)]ds = - \langle \bc_i\cdot\nabla_{\bx_i}\mu_t^{(N)},\Phi\rangle,\\
&\int_0^T \mathbb{E}[ \nabla_{\bc_i}\Phi(s,\bZ_s) \cdot \bH(s,\bZ_s)]ds =
- \left\langle\nabla_{\bc_i} \cdot \left[ \left(\G(.,\bz_i)+\frac{1}{N}\sum_{j=1}^N \F(\bz_i-\bz_j)\right) \mu_t^{(N)} \right],\Phi\right\rangle,\\
&\int_0^T \mathbb{E}[ \Delta_{\bc_i}\Phi(s,\bZ_s)]ds = \langle\Delta_{\bc_i}\mu_t^{(N)},\Phi\rangle.
\end{align*}
Since $\Phi$ does not have a compact support in time, integration by part requires to keep the boundary terms and the time partial derivative of $\mu_t^{N}$ defines as:
\begin{align*}
\int_0^T \mathbb{E}\left[ \frac{\partial}{\partial t}\Phi(s,\bZ_s)\right]ds = - \left\langle\partial_t \mu_t^{N}, \Phi\right\rangle + \mathbb{E}[ \Phi(t,\bZ_t)-\Phi(0,\bZ_0)].
\end{align*}
To sum up, thanks to the Itô's formula, we have obtained a weak form of the equation followed by the law $\mu_t^{(N)}$ of $\bZ(t)$:
\beq{eq:NJointWeak}
\resizebox{0.92\tw}{!}
{$ \ds
  \left\langle\partial_t \mu_t^{(N)} 
  + 
  \sum_{i=1}^N \bc_i\cdot\nabla_{\bx_i}\mu_t^{(N)} 
  + 
  \sum_{i=1}^N \nabla_{\bc_i}\cdot\left( \left[\G(.,\bz_i)+\frac{1}{N}\sum_{j=1}^N \F(\bz_i-\bz_j)\right] \mu_t^{(N)} \right),\Phi\right\rangle 
  = 
  \sigma ^2\left\langle\sum_{i=1}^N \Delta_{\bc_i}\mu_t^{(N)},\Phi\right\rangle.
$}
\eeq
So now, we have generalized the results given by Bolley in \cite{Bolley10} to a time dependent transport term $\G$.

However, equation \eqref{eq:NJointWeak} is a weak formulation of a PDE on the $N$-particles joint law, when what we are looking for is the equation ruling the one-particle law $\mu_t^{(1)}$, which is the marginal of $\mu_t^{(N)}$ for 
particle $1$. By integrating Eq.~\eqref{eq:NJointWeak} over all the particles but the first one, we get that, in the weak sense,  
$\mu_t^{(1)}$ follows: 
\beq{eq:FirstMarginal}
    \partial_t \mu_t^{(1)} + \bc \nabla_{\bx}\mu_t^{(1)} 
    +
    \nabla_{\bc}  \cdot \left( \G_t\mu_t^{(1)} + \int_{\bz_2 \in \mathbb{R}^{2d}} \!\!\!\!\!\!\!\!\!\!\!\!\!\! \F(\bz-\bz_2)\mu_t^{(2)}(\bz,\bz_2)  \right)
    \! = 
    \sigma ^2\Delta_{\bc}\mu_t^{(1)}.
\eeq
In this expression, $\mu_t^{(2)}$ is the 2-particles joint probability. 
In order to close equation \eqref{eq:FirstMarginal}, we would like to express it as a function of $\mu_t^{(1)}$.

To do so, we suppose that the initial data $\bZ_0$ are indistinguishable and all follow the same law $f_0$ on $\R^{2d}$. 
Then, we introduce an intermediate law, as the solution for $t>0$ and $(\bx,\bc)\in \mathbb{R}^{2d}$, of the following equation with initial data $f_0$:
\begin{equation}\label{eq_ft}
\frac{\partial}{\partial t} f_t + \bc \cdot \nabla_{\bx} f_t + \nabla_{\bc}\cdot \left[\left(\G + \F*f_t\right)f_t\right] = \sigma^2 \Delta_{\bc} f_t,
\end{equation}
where $\F$ and $\G$ are now supposed to be  Lipschitz functions with respect to the variable $\bx \in \mathbb{R}^{2d}$ and $\G$ is continuous in the time variable. 
Next, for $i \in 1,\dots,N,$ let $\bar \bz_i(t)$ be the solution of the following system with initial data $\bar \bz_i(0)=\bz_i(0)$:
\begin{equation}\label{syst_ft}
  \left\{
  \begin{aligned}
    &d \bar \bx_i(t)=\bar \bc_i(t) dt, \\
    & d \bar \bc_i(t)= \G(t,\bar \bz_i(t)) dt + \F*f_t(\bar \bz_i(t))dt+ \sqrt{2}\sigma d \bW_i(t).
   \end{aligned}
   \right.
\end{equation}
Then, the fictive particles $\bar \bz_i$ evolve in the field $\F*f_t$ generated by the distribution $f_t$, 
while the $\bz_i$ particles evolve in the $\F* \mu_t^N$ field, generated by the empirical measure $\mu_t^N$. 
Itô's formula gives once more the PDE followed by $\bar\bz_i$ in the weak sense, and we now wish to
show that this measure converges to $f_t$ when the number $N$ of particles tends to infinity. 

We denote $|(\bx,\bc)|=\sqrt{|\bx|^2+|\bc|^2}$ and for $p>1$ we define
$$\mathbb{P}_p(\mathbb{R}^{2d})=\left\{ \mu \text{ borelian probabilistic measures on } \mathbb{R}^{2d} \text{ such that } p\text{-momentum }  \int_{\mathbb{R}^{2d}} |(\bx,\bc)|^p d\mu(\bx,\bc) < \infty\right\}.$$
The Wasserstein distance of order $p$ between two measures $\mu$ and $\bar \mu$ of  $\mathbb{P}_p(\mathbb{R}^{2d})$  is defined by 
$$W_p(\mu,\bar \mu)=\inf_{\bZ,\bar \bZ} \sqrt[p]{\mathbb{E}\left[|\bZ-\bar \bZ|^p \right]},$$
where $\bZ$ and $\bar \bZ$ are stochastic variables of law  $\mu$ and $\bar \mu$ respectively. 
Then, following the lines of \cite{Bolley10},
\begin{thrm}\label{thrm-julie}
we obtain the explicit convergence rates: 
\begin{itemize}
\item[1)] $\ds W_2(\mu_t^{(1)},f_t)^2 \leqslant \mathbb{E}\left[\left|\bz_1(t)-\bar \bz_1(t)\right|^2\right] \leqslant \frac{C}{N}$,
\item[2)] $\ds W_2(\mu_t^{(k)},f_t^{\otimes k})^2 \leqslant \mathbb{E}\left[\left|\left(\bz_1(t),...,\bz_k(t)\right)-\left(\bar \bz_1(t),...,\bar \bz_k(t)\right)\right|^2\right] \leqslant \frac{Ck}{N}$,
\item[3)] Let $\Phi$ be a Lipschitz function in the second variable, then 
$$\mathbb{E}\left[\left| \int_{\mathbb{R}^{2d}} \Phi \mu_t^N - \int_{\mathbb{R}^{2d}} \Phi f_t\right| \right] \leqslant \frac{C}{N} ||\Phi||_1^2.$$
\end{itemize}
\end{thrm}
In other words, this means that:
\begin{itemize}
\item[1)] The one-particle law $\mu^{(1)}_t$ converges to $f_t$ in the Wasserstein distance when $N\rightarrow\infty$,
\item[2)] At the limit of an infinite number of particles, the chaos propagates in time; the particles remain uncorrelated during the whole dynamics: $\mu_t^{(k)}=f_t^{\otimes k}=\underbrace{f_t\otimes\dots\otimes f_t}_{\text{k times}}$. 
   In particular, one recovers the famous \textit{molecular chaos} assumption of Boltzmann:
   \beq{eq:MolecularChaos}
     \mu_t^{(2)} = \mu_t^{(1)}\otimes\mu_t^{(1)}.
   \eeq
\item[3)]The weak convergence of the empirical measure $\mu^N_t$ to  $f_t$.
\end{itemize}

Finally, equation \eqref{eq:FirstMarginal} is now closed rigorously thanks to the \textit{molecular chaos} propagation 
in the context of Lipschitz-regular interactions (external $\G$ or between particles $\F$), 
\cite{Villani02, Bolley10}. 
However, when the interactions are less regular, which is the case for the Boltzmann equation \eqref{eq:Boltzmann} below, 
an increasing number of positive results let us think that equation \eqref{eq:MolecularChaos} remains correct, \cite{Villani01,Lanford75}. Nonetheless, no rigorous demonstration is nowadays available.  
%\eitem 

\section{A population of particles in a turbulent fluid}
\label{sec:Macro}
In the previous section, a general kinetic equation has been derived for a population of ``particles'' (molecules, droplets, solid particles). As this point, one can be interested in deriving a two-way coupled system of kinetic equations for the carrying fluid and the particles. However, in \cite{doisneau_th13}, it has been shown in the context of nano-particles that such a derivation cannot be performed. Instead, we use the classical strategy of first deriving macroscopic equations for the fluid, and then coupling them to the particle equations, either microscopic or mesoscopic.
In the following, we first present the Euler and Navier-Stokes equations that can describe a carrying fluid, with an emphasis on the underlying assumptions at the kinetic level. In a context where dealing with the whole range of scales of the fluid is not accessible, we detail a general strategy for generating large-scale reduced-order models, and we show how it can be taken into account for the description of the particle dynamics at the microscopic level.

\subsection{Classical theories for macroscopic equations for the fluid}
In the context of gaz dynamics, in the limit of an infinite number of particles and when ignoring the stochastic 
subscale Brownian perturbations for the moment, equation \eqref{eq:FirstMarginal} becomes the Boltzmann equation:
\beq{eq:Boltzmann}
  \partial_t f + \bc\cdot\partial_{\bx} f + \partial_{\bc} \left(\Fext f\right) = \frac{1}{\Kn} Q(f,f),
\eeq
where $\Kn=\frac{\lambda}{L}$ is the Knudsen number, ratio between the mean free path $\lambda$ 
and a characteristic size of observation $L$, and where the quadratic collision operator $Q$ writes:
\beq{eq:OperateurCollision}
  Q(f,f_{*}) = \int_{\R^3_{\bc^*}} \int_{S^2_{\bn}} \left(f'f_*' - ff_*\right) 
  \aire{(\bc-\bc_*)\cdot\bn}\sigma(\aire{\bc-\bc_*},\bn) \; d\bn\; d\bc_*.
\eeq

\begin{rmrk}
  When considering a non self-interacting population of particles, its repartition function 
  also follows an equation of the \eqref{eq:Boltzmann} type, where the Knudsen number is infinite: 
  $\Kn = +\infty$. 
\end{rmrk}

\subsubsection{Euler equations}

For any PDF $f$, one can define its \textit{microscopic entropy} by $h = f\log f$. It can be understood as a local uncertainty rate.
Then, the \textit{macroscopic entropy} reads: 
$\H(t,\bx) = \int_{\R^3} f(t,\bx,\bc) \log f(t,\bx,\bc) \; d\bc,$ 
and one can show that when $f$ is a solution of the Boltzmann equation \eqref{eq:Boltzmann}, its macroscopic entropy 
decreases: $\dt{\H} \leq 0.$
When the minimum $\H_{\text{min}}$ is reached, $\log f$ must be a \textit{collision invariant} and this implies that
the velocity distribution $f$ is a Maxwellian: 
\beq{eq:DistributionEquilibre}
  \feq(\bc) = \exp\left(a_0 + \boldsymbol{a_1}\cdot\bc + a_2\frac{|\bc^2|}{2}\right)
         = f_0 \exp\left(-\frac{|\bc-\bu|^2}{2\beta}\right).
\eeq
The Maxwellian distribution being perfectly defined by its three first moments 
$m_k = \int_{\R_{\bc}^3} \bc^k f(\bc)d\bc$, $k=0,1,2$, the evolution of the Boltzmann equation 
\eqref{eq:Boltzmann} at isentropic thermodynamic equilibrium $\H = \H_{\text{min}}$ is given by the system of its three first moments, 
which closes into the Euler equations: 
\beq{eq:Euler}
\left\{
  \begin{array}{cccl}
    \dd{\rho}{t}     & + & \Div{\rho \bu}                 & = 0,\\
    \\
    \dd{\rho \bu}{t} & + & \Div{\rho \bu\otimes\bu + p\I} & = 0, \\
    \\
    \dd{\rho E}{t}        & + & \Div{(\rho E+p)\bu}       & = 0.
  \end{array}
\right.
\eeq

\subsubsection{Navier-Stokes-Fourier equations}

In the previous paragraph, we have clearly stated that collisions occur everywhere at all time, or, to reformulate, that the Knudsen number remains null: $\Kn=0$. 
In reality, it is often very small but strictly positive. Then, we look at near equilibrium regimes by stating $\eps = \Kn$ and looking for 
an expansion of $f$ in $\eps$: the Chapman-Enskog expansion. 
At first order, $f = f_0 + \eps f_1 + \circ(\eps)$, 
which, at orders $1/\eps$ and $1$, gives in \eqref{eq:Boltzmann}:
$$
  Q(f_0,f_0) = 0 \quad \Longrightarrow \quad f_0 = \feq 
  \quad \quad \text{ and } \quad \quad
  \partial_t \feq + \bc\cdot\partial_{\bx}\feq = Q(\feq,f_1) + Q(f_1,\feq).
$$
The latest equation is an integral equation in $f_1$ which might be completely solved. For a monoatomic gas of atoms of mass $m$ and radius $r$, 
the three first moments of $f$ verify the following Navier-Stokes equations 
\cite{ChapmanCowling}:

\begin{minipage}{0.4\tw}
\beq{eq:NS}
\resizebox{0.85\tw}{!}
{$ %footnotesize
  \left\{
    \begin{array}{cccc}
    \ds \dd{\rho}{t}      &+&  \Div{\rho \bu}                      &= 0,\\
    \\
    \ds \dd{\rho \bu}{t}  &+&  \Div{\rho \bu\otimes\bu + \P}       &= 0, \\
    \\
    \ds \dd{\rho E}{t}    &+&  \Div{\rho E \bu + \P\cdot\bu + \bq} &= 0.
    \end{array}
  \right.
$}
\eeq
\end{minipage}
\begin{minipage}{0.02\tw}
~
\end{minipage}
\begin{minipage}{0.05\tw}
%\footnotesize
\vspace{1.5em}
where
\end{minipage}
\begin{minipage}{0.02\tw}
~
\end{minipage}
\begin{minipage}{0.44\tw}
\vspace{0.7em}
\begin{equation*}
\resizebox{0.95\tw}{!}{$ %\footnotesize
 \left\{
   \begin{array}{ll}
     \ds \P &= \ds p\I - \frac{\mu}{2}\left(\grad\bu + \grad\bu^t-\frac{2}{3}\div\bu \I\right),\\
      \\
     \bq &= \eps\int_{\R^3_{\bc}} m \frac{(c-u)^2}{2} (\bc-\bu) f_1 d\bc = -\lambda\grad T,
   \end{array}
 \right.
$}
\end{equation*}
\end{minipage}\\
\\
and $\mu$ and $\lambda$ are respectively the \textit{viscosity} and the \textit{thermal conductivity}. 

\begin{rmrk}
This last system can be obtained rigorously from the Boltzmann equation \eqref{eq:Boltzmann} 
in the restrictive context of monoatomic gases with 
$\mu = \frac{5m}{6r^2}\sqrt{\frac{k_B T}{\pi m}}$ and 
$\lambda = \frac{225}{512 r^2} \sqrt{\frac{m k_B T}{\pi}}$. 
However, a similar system of PDEs can be obtained by considering the 
conservative principles of mass, momentum and total energy, 
added with constitutive equations of the considered fluid, which provide heuristic 
laws of the viscosity $\mu$ and the thermal conductivity $\lambda$. 
\end{rmrk}

\subsection{Large-scale reduced-order models}
\label{lsrom}
\subsubsection{Properties of turbulence}\label{sssec:Turbulence}

Turbulence is a particular type of flows which can not be rigorously defined. 
The easiest way to define it is by using the metric of the Reynolds number: 
$\text{Re}=\frac{u_{f,0}L_0}{\nu}$,
where \(u_{f,0}\) is a characteristic speed of the fluid, \(L_0\) is a characteristic length scale of the system and 
$\nu=\frac{\mu}{\rho}$ is the kinematic viscosity of the fluid. 
We will say that a fluid exhibits a \textit{turbulent behavior}, when its Reynolds number is high. 
The limit Reynolds number depends on the considered experiment and on the operating condition. 
However, the flow is generally turbulent when $\text{Re} >> 10^3$. 

Turbulent flows share in common their chaotic behavior. 
For deterministic systems, there are multiple definitions of chaos, but in this context we choose to say that turbulent flows all are~:
\begin{itemize}
\item \textit{highly sensitive to the initial conditions of the system}. 
  The present determines the future, but the approximate present does not approximately describe the future. 
  For instance, we say that $x_{0}$ is a highly sensitive initial conditions, if for all $L_0 > M >0$ and for all $\delta>0$, 
  there exists another close initial data $y_{0}$ and an arbitrary time $t>0$ such that 
\[
|x_{0}-y_{0}| < \delta \text{ and  } |x(t;x_{0})-y(t;y_{0})| \geq M.
\]
\item \textit{topologically transitive}, 
  in the sense that for every pair of non-empty open sets $U\subset X$ and $V\subset X$, 
  there is an arbitrary time $t>0$ such that
\[
\{ x(t;x_{0}) \in X : x_{0} \in U \} \cap V \neq \emptyset.
\]
\end{itemize}

From an experimental point of view, some observations have been made on turbulent fluid flows.
The main ones are expressed by Kolmogorov \cite[p.190]{Pope01}.
\begin{itemize}
\item At sufficiently high Reynolds number, the small-scale turbulent motions are statistically isotropic. 
  They follow a universal form that is uniquely determined by the viscosity \(\nu\) and the energy dissipation \(\varepsilon\).
\item The viscosity also defines a cut-off size $\eta_K$, called the \textit{Kolmogorov scale}, below which all the inertia of the flow is dissipated. 
\item Between the characteristic length $L_0$ and $\eta_K$, there is an intermediate range of scales, 
  called \textit{the inertial range}, where the statistics of motion have a universal form 
  that is uniquely determined by the dissipation \(\varepsilon\) and is independent of the viscosity \(\nu\). 
  Through dimensional analysis, we get that within this 
   range, the \textit{turbulent kinetic energy} decreases as: 
  \(E\left(\left|\boldsymbol{k}\right|\right)\propto\left|\boldsymbol{k}\right|^{-\frac{5}{3}}\), with \(\boldsymbol{k}\) the wavenumber.
\end{itemize}
%o[[Est-ce qu'on detaille ?]]\\
%o[[ergodicit\'e, limitation des petites \'echelles par la viscosit\'e]]

\subsubsection{Reduced description of turbulence}
\label{ssec:ProbaMapPope}

It is commonly admitted that the macroscopic Navier-Stokes equations contain the turbulence defined above, in the sense that 
these equations present solutions which have all the properties listed in paragraph \ref{sssec:Turbulence}.
Nonetheless, in practice the domain size, denoted by $|\mathbb{X}|$, and the dissipative cut-off scale $\eta_K$, may be separated by many orders of magnitude. 
In this context, the Direct Numerical Simulation of the Navier-Stokes equations is rapidly unreachable, since the number of needed computational cells 
will be at least of the order of $(\eta_K/|\mathbb{X}|)^3$, not speaking about the generally necessary high number of degrees of freedom per cell. 

Therefore, while staying very generic, 
we consider a decomposition of the solution into a \textit{significant part} and a \textit{residual}:
if $\phi$ is a quantity of interest, we consider its reduction $\moy{\phi}$ on the space of significant data
and thus write $ \phi = \moy{\phi} + \phi'$.
This significant part could be an \textit{ensemble average}, a \textit{filtering}, 
a spatial or a temporal average or even a modal decomposition. 
The goal is always to reduce the size of the information needed to entirely represent the chosen significant part, 
hence the name \textit{reduced-order} model. 

Now, the reduction operator $\moy{\ \cdot\ }$ is applied directly 
on the macroscopic equations Eq.~\eqref{eq:Euler}-\eqref{eq:NS}. 
For example, when considering the incompressible version of the Navier-Stokes equation, assuming commutativity between all implied linear operators, one gets: 
\begin{equation}
\begin{array}{rcl}
\nabla\cdot\overline{\boldsymbol{u}_f} & = & 0,\\
\partial_t\left(\overline{\boldsymbol{u}_f}\right)+\overline{\left(\boldsymbol{u}_f\cdot\nabla\right)\boldsymbol{u}_f} - \nu\nabla^2\overline{\boldsymbol{u}_f} & = & -\frac{1}{\rho_f} \nabla \overline{p},
 \end{array}
\label{eq:NSconv}
\end{equation}
with \(\boldsymbol{u}_f\) the fluid velocity, \(\rho_f\) its density (constant for incompressible fluids), 
\(\nu\) its kinematic viscosity and \(p\) the pressure field.

The main difficulty now lies in the reduction of the non-linear terms. 
Indeed, nothing indicates that there exists an application giving 
\(\overline{\left(\boldsymbol{u}_f\cdot\nabla\right)\boldsymbol{u}_f}\) as a function of \(\overline{\boldsymbol{u}_f}\).
Thus, Eq.~\eqref{eq:NSconv} is not meaningful in term of the significant unknown \(\overline{\boldsymbol{u}_f}\).
To overcome this difficulty, the main idea is to define a more complex application which gives
multiple possibilities to the relation between 
\(\overline{\left(\boldsymbol{u}_f\cdot\nabla\right)\boldsymbol{u}_f}\) and \(\overline{\boldsymbol{u}_f}\).
This is done by adding a hidden variable $\omega$, which encodes all the complexity of \(\overline{\left(\boldsymbol{u}_f\cdot\nabla\right)\boldsymbol{u}_f}\) inside an application $\mathcal{F}$ and a space of possibilities $\Omega$ 
in the following way
\beq{eq:F}
\mathcal{F}:
 \left\{
   \begin{array}{rcl}
     \Omega \times \R^3_{\bu} &\rightarrow & \R^3_{\bu}\\
     (\omega, \overline{\boldsymbol{u}_f}) &\mapsto& \overline{\left(\boldsymbol{u}_f\cdot\nabla\right)\boldsymbol{u}_f} (\omega).
   \end{array}
 \right.
\eeq
Of course, the definition of $\mathcal{F}$ strongly depends on the choice of the reduction operator $\moy{\ \cdot\ }$.
Next, an elegant way to move forward is now to define $\Omega$ as a probability space, see \cite{pope2010self}. 
Then, two main techniques emerge~:
\begin{itemize}
\item by drawing many particular $\omega$, thus giving a random modeling of the unknown term   \(\overline{\left(\boldsymbol{u}_f\cdot\nabla\right)\boldsymbol{u}_f}\) through $\mathcal{F}$, compute many trajectories of the process \(\overline{\boldsymbol{u}_f}\),

\item considering the statistics or moments of the random variable $\mathcal{F}$, and solve for the evolution of the moments of the random variable \(\overline{\boldsymbol{u}_f}\).
\end{itemize}

The advantage of the first approach is to preserve the properties of a trajectory of the process 
\(\overline{\boldsymbol{u}_f}\), which is still the solution of a PDE. Thereby, 
the random variable \(\overline{\boldsymbol{u}_f}\) lies in a large dimensional probability space, which requires a very large number
of such succession of draws to hope for some meaningful statistics. 
On the contrary, solving for the evolution of the means of \(\overline{\boldsymbol{u}_f}\) 
does not preserve the trajectories of the process, but it gives correct estimators and  statistics on the 
general behavior of the gaseous velocity field.

\subsubsection{Closures}

The obtained reduced-order system as in Eq.~\eqref{eq:NSconv} is closed by making a calculable choice on $\mathcal{F}$. Three strategies can be found in the literature for this choice, as depicted in \cite{sagaut2006large, Pope01}:
\begin{itemize}
\item the functional approach: starting from the fact that the regularized version of the flow field will dissipate less energy than the real turbulent flow field does, the unresolved scales can be modeled  in a first approximation by an additional diffusion process, consistently with the theory of turbulence described in paragraph \ref{sssec:Turbulence}:
$$
\overline{\left(\boldsymbol{u}_f\cdot\nabla\right)\boldsymbol{u}_f} - \left(\uf\cdot\nabla\right)\uf
\approx
-\mu^{turb} \nabla^2 \uf.
$$
Here, $\mu^{turb}$ is an additional \textit{turbulent viscosity}. In the case of filtering procedures, this viscosity depends on the filter size such that it vanishes for full-resolution \cite{smagorinsky1963general,nicoud1999}.  As such models can depend on empirical constants, dynamic procedures were also proposed to get the better estimate of theses constants (see \cite{germano1991dynamic}).

\item The structural approach: instead of simply recovering a global property of the unresolved information, structural methods aim at capturing the SGS tensor structure (see \cite{Bardina:1980}).

\item the "pragmatic" approach: starting from the idea that it is hard to distinguish unresolved scales effects from numerical dissipation, some authors propose to integrate effects of unresolved scales through  the numerical schemes (see \cite{grinstein:2002}).
\end{itemize}
%\section{Inertial particles in turbulent flow}
%\label{sec:TurbulentParticles}

\subsection{Particles in turbulence}

\subsubsection{Reduced LES models}

The fluid velocity at the location of the particle appears in the expression of the particle acceleration modeled by Stokes drag law:
\beq{eq:Stokes}
  d\bC_i = \dfrac{\bu_f(t,\bX_i)-\bC_i}{\tau_p} dt, 
\eeq
\(i\in [\![ 1,N ]\!]\), and $\tau_p$ being a characteristic relaxation time of the particule toward the underlying velocity field.
However, in every LES model existing up to now, only a regularized version of the fluid velocity is computed. Thus, a closure on the fluid velocity seen by the particle is required in order to provide a consistant LES model for the disperse phase. Ideally, this model has to be in agreement with the probability space of the random variable \(\mathcal{F}\) seen by the inertial particles on the fluid flow.

Up to now, very similarly to the models developed for the fluid flow, the main strategies have been to compensate second order moments of the the particle density distributions by the adjunction of energy in the form of Wiener processes (see \cite{bini2007particle, fede2006numerical, minier2004pdf, reeks1977dispersion, shotorban2007eulerian, shotorban2006stochastic}). 
In its general from, this can be represented by the stochastic differential equation \eqref{eq:Ito}:

\begin{equation}\label{eq:Ito}
d\bZ_t=\mu_t \bZ_t dt + \sigma_td\bW_t,
\end{equation}
with \(\bW_t\) a Wiener process, \(\bZ_t\) the state vector of the particle, \(\mu_t\) the drift and \(\sigma_t\) the diffusion coefficient. It is to be noted that in most models, the Wiener process only acts on one variable of the particle~: either its position, or its velocity, or an other intermediate variable like the velocity seen by the particle. The next section shows that in the context of equation \eqref{eq:Ito}, where the closure has been chosen in the form of a Wiener process, the derivation of a mesoscopic equation for the disperse phase is not a major difficulty.

\section{Consistency of modeling approaches with numerical cases}
\label{sec:Numerics}

Sections \ref{sec:MicroMeso} and \ref{sec:Macro} were mainly focused on providing a meaningful formalism for reduced multiphase flow simulations in agreement with mathematical consistency and physical literature. In this context, we conclude that an appropriate formalism to describe a fluid in a Large-scale reduced order in section \ref{lsrom} is the self-conditioned structure proposed by \cite{pope2010self} and formalized Eq.~\eqref{eq:F}. In a nutshell, the evolution of the large scale of the flow must be obtained as the expectation of all possible unresolved scales of the flow compatible with the resolved large scales. 

Applying this formalism with the full resolution of Navier-Stokes is not easy because it is not straightforward to control large scales and unresolved scales separately. An interesting alternative that has been widely used in the literature is to rely on synthetic turbulence: by means of a summation of analytic modes, and under the constraint of specific spectral distribution and representation, one can expect to reproduce the main characteristics of the turbulence, even without verifying Navier-Stokes equations.
In this section, we investigate the use of such analytic representation from 1D to 3D, and we show what is the minimal representation that can be envisaged.

% Thus, we will start by working on a synthetic fluid flow reproducing some of the essential properties of the solutions of Navier-Stokes equations (such as incompressibility) and having a clear and simple self-conditioned structure. The design of this fluid flow and the characterisation of some of its properties will be detailed in this section.
% This section aims at building a numerical framework, induced from the previous results, in order to numerically validate  models according to this reduction formalism. The objective is to assess the relation between reduced fluid flow descriptions and corresponding disperse phase closure models. 
%From what is known of turbulence today, it seems difficult to start by working directly on the numerical solutions of Navier-Stokes equations if one wants to investigate the impact of a self-conditioned field framework on the disperse phase modeling. Thus, we will start by working on a synthetic fluid flow reproducing some of the essential properties of the solutions of Navier-Stokes equations (such as incompressibility) and having a clear and simple self-conditioned structure. The design of this fluid flow and the characterisation of some of its properties will be detailed in this section.

\subsection{Synthetic turbulence}
The synthetic flow field has been designed in order to reproduce somehow the dynamics that could be expected from a self-conditioned LES flow field simulation (\hspace{1sp}\cite{ijzermans2010segregation, kraichnan1970diffusion}).
It is represented by a sparse matrix of spectral modes (Eq.~\eqref{eq:fluidModes}) chosen according to the energy density given by Pope's spectrum in Eq.~\eqref{eq:PopeSpectrum} (see \cite[p.232]{Pope01}) with Eq.~\eqref{eq:modesChoice}.

\begin{equation}
\boldsymbol{u}_f\left(t,\boldsymbol{z}\right)=\sum_{n=0}^{N}\boldsymbol{a}_n\cos\left(\omega_n t+\boldsymbol{k}_n\cdot\boldsymbol{z}+\varphi_n\right)
\label{eq:fluidModes}
\end{equation}

\begin{equation}
E\left(\left|\boldsymbol{k}\right|\right)=\frac{9}{4}\frac{\epsilon^{\frac{2}{3}}}{\left|\boldsymbol{k}\right|^{\frac{5}{3}}}
\left(\frac{\left|\boldsymbol{k}\right|/k_0}{\left[\left(\left|\boldsymbol{k}\right|/k_0\right)^2+6.78\right]^{1/2}}\right)^{\frac{11}{3}}
\exp\left(-5.2\left(\sqrt[4]{\left(\left|\boldsymbol{k}\right|\eta\right)^4+0.4^4}-0.4\right)\right)
\label{eq:PopeSpectrum}
\end{equation}

\begin{equation}
\int_0^{\left|\boldsymbol{k}\right|_n}E_{3D}(\left|\boldsymbol{k}\right|)d\left|\boldsymbol{k}\right|=\frac{3}{2}u_0^2\frac{\left(2n-1\right)}{2N}.
\label{eq:modesChoice}
\end{equation}

The amplitude of the modes is chosen according to the distribution \(\left|\boldsymbol{a}_n\right|\sim\mathcal{N}\left(0,\frac{2u_0^2}{N}\right)\).

Following \cite{ijzermans2010segregation}, the spectral components of the energy spectrum are chosen in order to respect the numerical simulations performed in \cite{hunt1987big}, which show that it seems sensible to approximate \(E\left(\left|\boldsymbol{k}\right|,\omega\right)\) by~:
\begin{equation}
E_{\omega}\left(\left|\boldsymbol{k}\right|\right)=\frac{E\left(\left|\boldsymbol{k}\right|\right)}{\sqrt{2\pi}\left(a\left|\boldsymbol{k}\right|u_0\right)}\exp\left(-\frac{\omega^2}{2\left(a\left|\boldsymbol{k}\right|u_0\right)^2}\right),
\label{eq:Hunt}
\end{equation}
with \(a\in\left[0.4,0.51\right]\) depending on the wavenumber and the integral length scale (see \cite{hunt1987big}).
For the numerical simulations, the random number generator chosen is \texttt{ran2} presented in \cite{press1989numerical}.
The numerical values are chosen such that \(a=0.5\), \(u_0=1\)~m.s\(^{-1}\) and \(k_0=1\)~m\(^{-1}\). 
The particle evolution is computed using Runge-Kutta scheme of order four.

The evolution of the particles on the fluid is computed by the linearised Stokes drag law in Eq.~\eqref{eq:Stokes}, with the expression of \(\boldsymbol{u}_f\) given in Eq.~\eqref{eq:fluidModes}.

\begin{figure}[ht]
\centering
\includegraphics[width=0.8\textwidth]{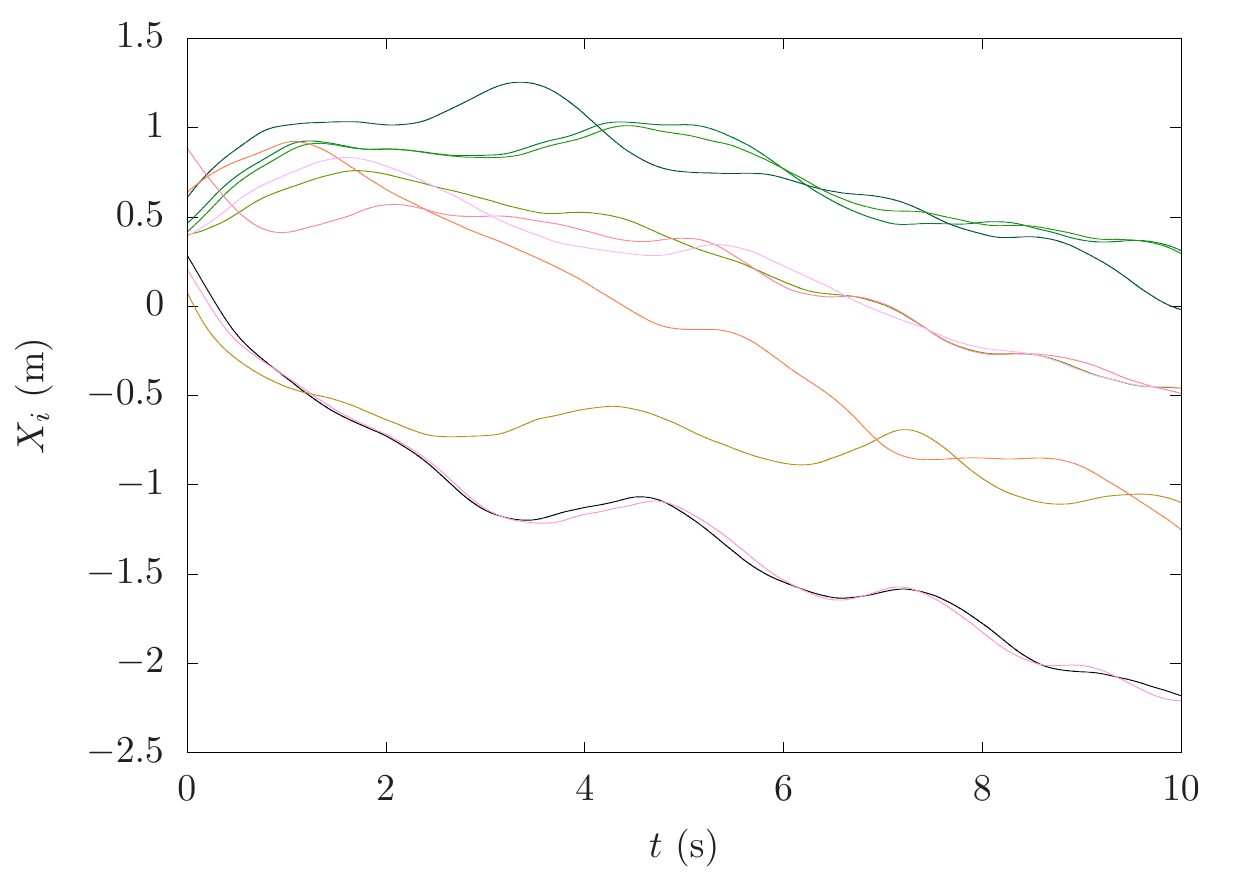}
\caption{Position $X_{i}$ in meters (m) of particles according to time in seconds (s) and colored by their number $i\in [\![0, 9]\!]$. Time step of \(dt=0.001s\) for particle time relaxation constant of \(\tau_p=1s\). }
\label{fig:part_xt-1D}
\end{figure}

For numerical simplicity, we first start by performing one-dimensional simulations. In one dimension, a realization of the evolution of the particles submitted to a random fluid is given in Fig.~\ref{fig:part_xt-1D}. 
Although the initial positons of the particles are random and uniformly distributed on a segment, their trajectories seem very limited. They look more like oscillations around a mean drift rather than dispersion.
Furthermore, when observing the evolution of the variance in a one-dimensional space for 10\(^{4}\) particles,
see Fig.~\ref{fig:stat-xpn-var-1D}, we see that it seems bounded for this case 
and that it is highly dependent on the underlying fluid fluctuations.

This kind of behavior is not consistent with the properties of turbulence and the expected behavior of particles 
in a turbulent flow: we would rather expect a dispersion behavior \textit{similar} to diffusion 
(see for instance \cite{richardson1926atmospheric}). 
Since the stochastic models of the literature have a first order effect on the second order moments of the measure 
of the disperse phase, it is essential to work on a numerical setup which preserves the basic properties of turbulent flows for realizations of the second order moments of the measure of the disperse phase. Hence, it is of prime importance to understand why such a behavior is observed on the simple fluid model we have chosen if we want to use it for reproducing and understanding the dynamic of inertial particles on fluids described by Navier-Stokes kind of equations.

\begin{figure}[!htp]
\centering
 \subfigure[1D]          
 	{\includegraphics[width=0.6\textwidth]{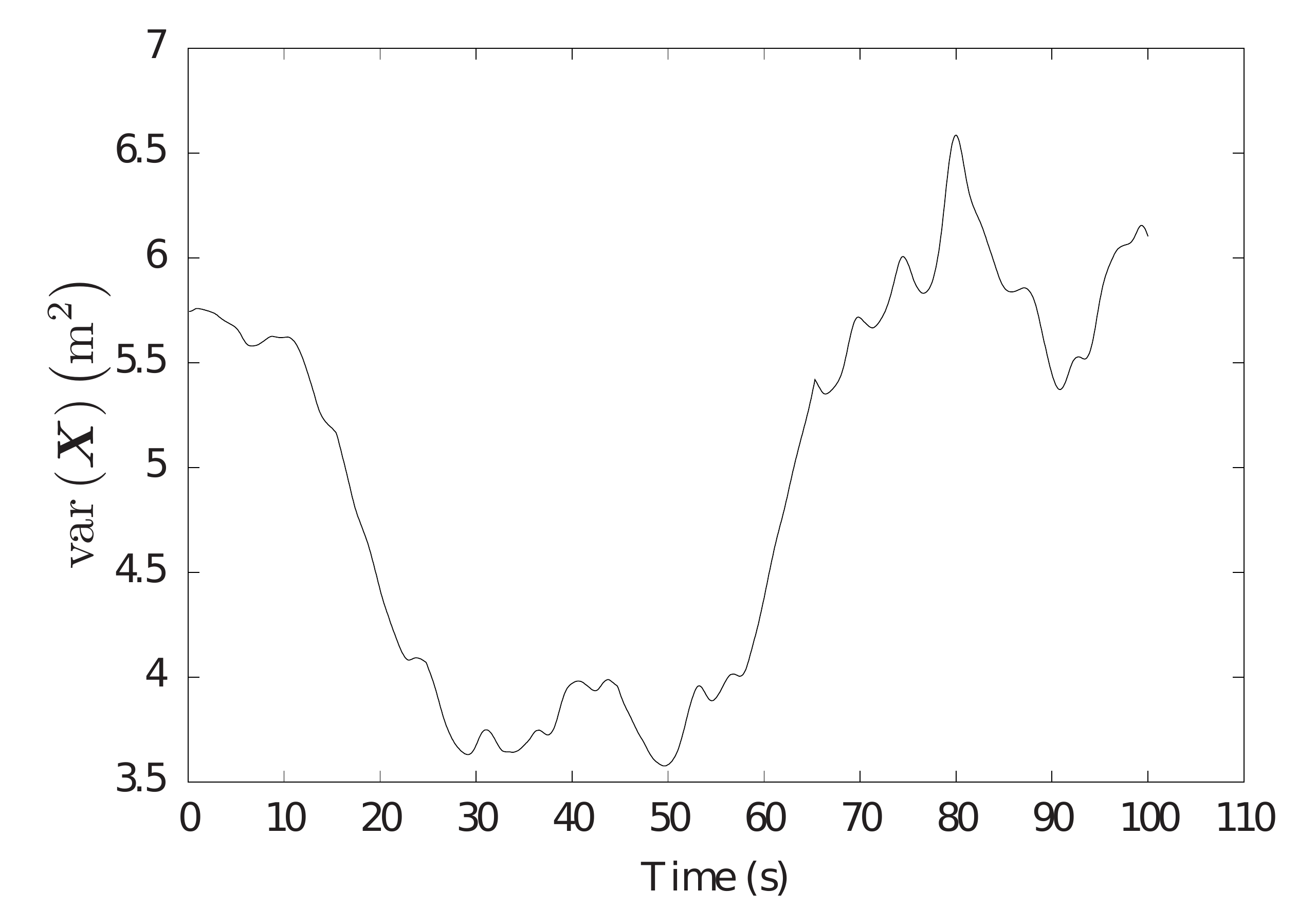}
	\label{fig:stat-xpn-var-1D} }
 \subfigure[2D]      
 	{\includegraphics[width=0.6\textwidth]{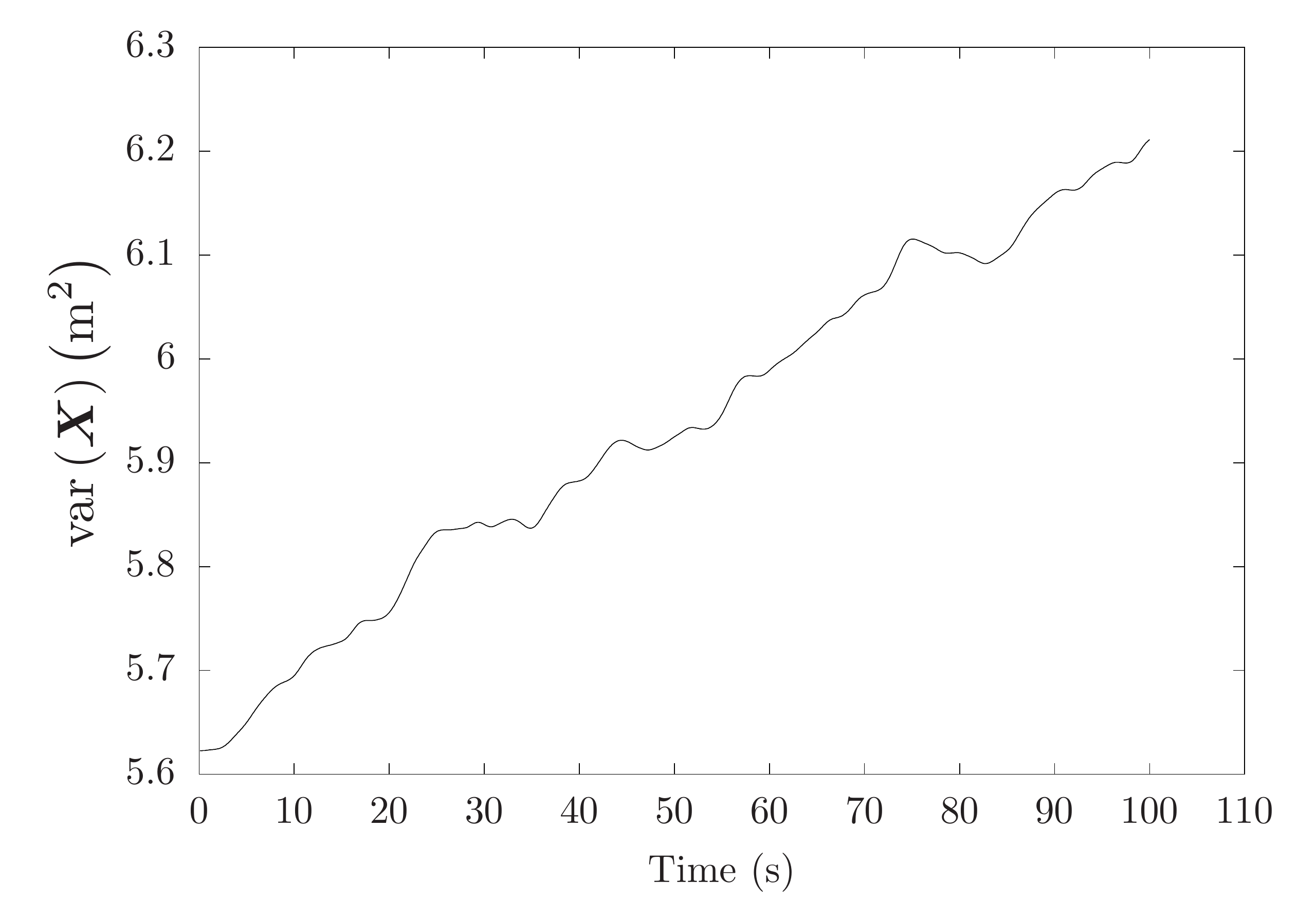}
	\label{fig:stat-xpn-var-2D}}
\caption{Evolution of the variance in position of 10\(^{4}\) particles of relaxation time constant of \(\tau_p=1\)~s in different dimensionalities on one fluid flow realization according to time (s).}
\end{figure}
\begin{figure}\ContinuedFloat
 \subfigure[3D]      
 	{\includegraphics[width=0.6\textwidth]{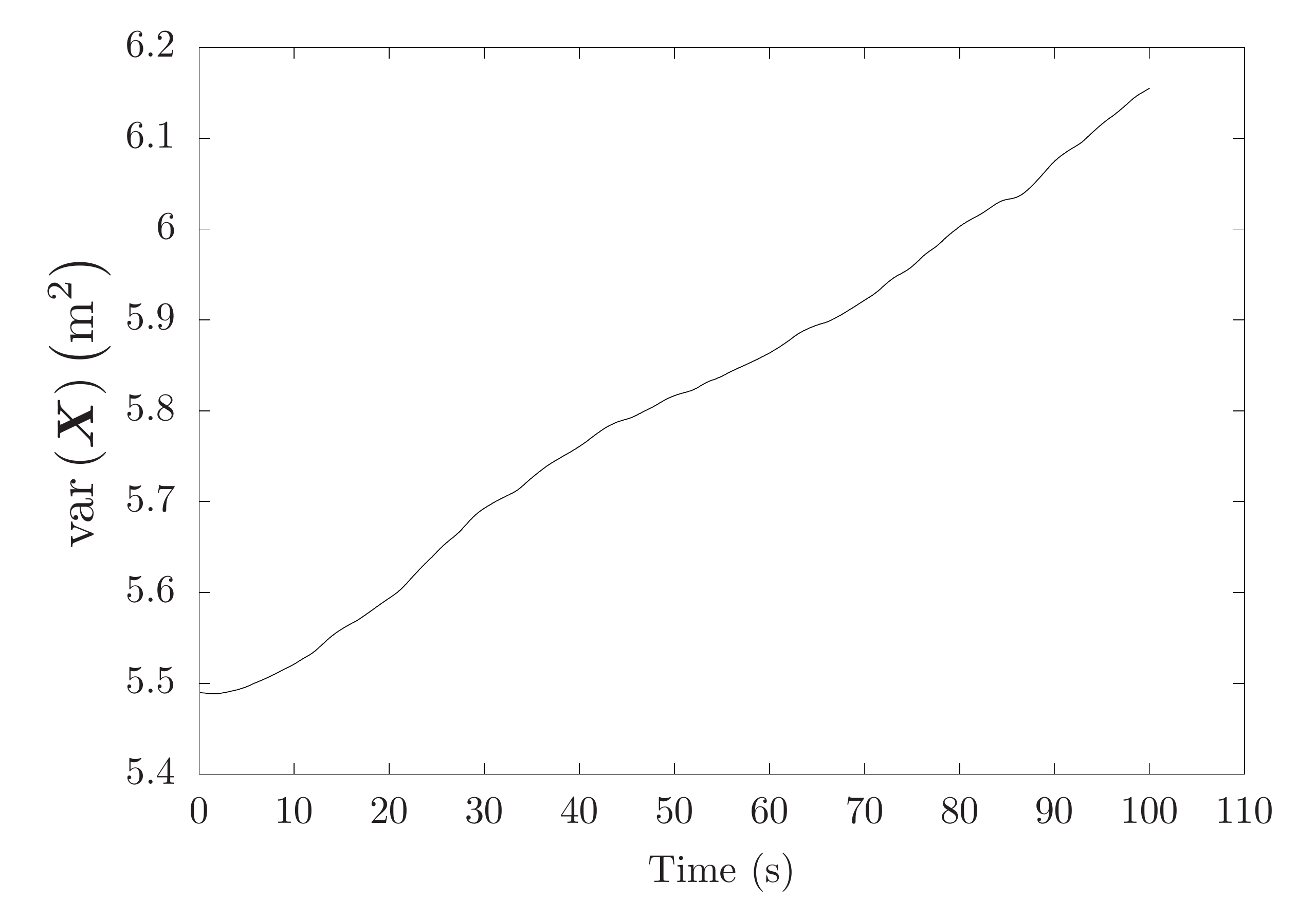}
	\label{fig:stat-xpn-var-3D}}
\caption{Evolution of the variance in position of 10\(^{4}\) particles of relaxation time constant of \(\tau_p=1\)~s in different dimensionalities on one fluid flow realization according to time (s).}
\label{dispersion-dim}
\end{figure}

\subsection{Simplified one-dimensional case}

As explained above, Fig.~\ref{fig:part_xt-1D} enlightens an unexpected behavior in one dimensional case. Let us start by looking if it is possible to understand this behavior on a simplified case where the fluid is only represented by one sine. We have the particle evolution in Eq.~\eqref{eq:equation1D} and the reduced evolution in Eq.~\eqref{eq:equation1Dreduced}.

\begin{equation}
 \left\{
\begin{array}{l}
dX_i=C_idt,\\
dC_i=\dfrac{a\sin\left(2\pi\left(\omega t + k X_i\right)+\varphi\right)-C_i}{\tau_p},
\end{array}
 \right.
 \quad i=1,\dots,N.
\label{eq:equation1D}
\end{equation}

\begin{equation}
 \left\{
\begin{array}{l}
dX_i'=C_i'dt,\\
dC_i'=\dfrac{a'\sin\left(2\pi X_i'\right)+\omega-C_i'}{\tau_p},
\end{array}
 \right.
  \quad i=1,\dots,N.
\label{eq:equation1Dreduced}
\end{equation}

We will prove the following result : 
\begin{prpstn} Particles under dynamics described in Eq.~\eqref{eq:equation1D} will follow an increasing signal, incompatible with an expected diffusive behavior.
\end{prpstn}

\begin{proof}
In order to prove this result, we can first study the system \eqref{eq:equation1Dreduced}. This 
system is autonomous in dimension 2, so by the Poincar\'e-Bendixon theorem, only three cases are possible~: 
\begin{itemize}
\item The trajectories are unbounded,
\item The trajectories converge to a point,
\item The trajectories converge to a limit cycle.
\end{itemize}

Let us now try to characterize these behaviors more precisely.

Let \(i\in[\![1,N]\!]\). Define $C_{\max} = |a'|+|\omega|+1$. Suppose that at time $t^{*}$ we have $C_{i}'(t^{*}) > C_{\max}$, then by continuity during a time $\delta$ we have $C_{i}'(t^{*}+t)\geq C_{\max}$ for all $t\in [0, \delta]$. Thus
\begin{eqnarray*}
C_{\max}\leq C_{i}(t^{*}+\delta)
&=& C_{i}'(t^{*}) + \int_{0}^{\delta} \frac{a'\sin\left(2\pi X_i'(t^{*}+s)\right)+\omega-C_i'(t^{*}+s)}{\tau_p}ds\\
&\leq& C_{i}'(t^{*}) + \int_{0}^{\delta} \frac{|a'|+|\omega|-C_{\max}}{\tau_p}ds\\
&\leq& C_{i}'(t^{*}) +\delta\frac{|a'|+|\omega|-C_{\max}}{\tau_p}
\end{eqnarray*}
or simply
%$v_{\max} (1+\frac{\delta}{\tau_{p}})\leq v_{p}'(t^{*}) + \delta\frac{|a'|+|\omega|}{\tau_p}$
$\delta \leq (C_{i}'(t^{*})-C_{\max})\tau_p$. It proves that in finite time, the solution falls under $C_{\max}$. Then we have proved that for all trajectories, there exists a time $t_{\max}$ where $C_{i}'(t_{\max})< C_{\max}$.

Denote $C_{\min}= -C_{\max}$ and use again the time $t^{*}$ with symmetric definition, we obtain
\begin{eqnarray*}
C_{\min}\geq C_{i}(t^{*}+\delta)
&=& C_{i}'(t^{*}) + \int_{0}^{\delta} \frac{a'\sin\left(2\pi X_i'(t^{*}+s)\right)+\omega-C_i'(t^{*}+s)}{\tau_p}ds\\
&\geq& C_{i}'(t^{*}) + \int_{0}^{\delta} \frac{-|a'|-|\omega|-C_{\min}}{\tau_p}ds\\
&\geq& C_{i}'(t^{*}) +\frac{\delta}{\tau_p}
\end{eqnarray*}
and $\delta \leq (C_{\min}-C_{i}'(t^{*}))\tau_p$. It proves that in finite time, the solution rises above $v_{\min}$. Then we have proved that for all trajectories, there exists a time $t_{\min}$ where $C_{i}'(t_{\min})> C_{\min}$.

Finally we can suppose that for all trajectories, the speed $C_{i}' \in [-C_{\max},C_{\max}]$ after some transitory time. In fact -with the same procedure- we can prove that $C_{i}' \in [\omega-|a'|,\omega+|a'|]$.

Thus if $\omega > |a'|$ then it proves that the speed stays strictly greater than $\omega-|a'|>\varepsilon>0$, and thus the trajectories cannot be bounded. In order to prove that the particles will follow an increasing signal, we have  to study the difference with this linear growing.

Denote $Y_{i}'(t) = X_{i}'(t) - \omega t$ then $V_{i}'(t) =C_{i}'(t)-\omega= \frac{dY_{i}'}{dt} \in [-|a'|,|a'|]$.

\begin{align}
dY_i'&=(C_i'-\omega)dt = V_{i}' dt,\\
dV_i'&=\frac{a'\sin\left(2\pi Y_i'+2\pi\omega t\right)-V_i'}{\tau_p}dt.
\end{align}

We can see that $V_{i}'$ cannot converge to a constant $\bar{V}$, because there is no solution to $\bar{V} = a' \sin(2\pi t (\bar{V}+\omega))$ (except $\omega = \bar{V} = 0$). Since $V_{i}'$ cannot converge to a constant while staying in a compact, it is non-monotonous. Denote $T_{+}$ a moment where 
$\frac{dV_{i}'}{dt}$ changes its sign (without loss of generality, suppose it changes from $>0$ to $<0$), i.e.
\[
a' \sin (2\pi (Y_{i}'(T_{+}) + \omega T_{+})) = V_{i}'(T_{+}) =: V_{+}.
\]
Thus 
\[
\tau_{p}\frac{d^{2}V_{i}'}{dt^{2}} (T_{+})=  2\pi a'\cos\left(2\pi Y_i'(T_{+})+2\pi\omega T_{+}\right)\left( V_{i}'(T_{+}) + \omega \right)-\frac{dV_i'}{dt}(T_{+})
\]
and in particular
\[
\dfrac{\dfrac{d^{2}V_{i}'}{dt^{2}} (T_{+})}{a'\cos\left(2\pi Y_i'(T_{+})+2\pi\omega T_{+}\right)}= \frac{2\pi}{\tau_{p}}\left( V_{+}+ \omega \right)
\]

Suppose $a'>0$ to simplify.

The quantity $\frac{dV_{i}'}{dt}$ changes from $>0$ to $<0$, thus the second derivative is negative, so at a given time, there is a local maximum, and during a period $[T_{+},T_{+}+T]$, $V_{i}'$ is decreasing and we have also $2\pi(Y_{i}'(T_{+})+\omega T_{+}) \in [\pi/2,3\pi/2] \mod 2\pi$. Or simply $2\pi(Y_{i}'(T_{+})+\omega T_{+}) = 2k\pi + \pi/2 + \varepsilon \pi$ with $\varepsilon \in [0,1]$.
% and we have 
%\begin{eqnarray*}
%0
%&\geq& \int_{T+}^{T} \frac{a'\sin\left(2\pi y_p'(t)+2\pi\omega t\right)-w_p'(t)}{\tau_p}dt\\
%&\geq& \int_{T+}^{T} \frac{a'\sin\left(2\pi y_p'(t)+2\pi\omega t\right)-W_{+}}{\tau_p}dt\\
%&\geq& \int_{T+}^{T} \frac{a'(\sin\left(2\pi y_p'(t)+2\pi\omega t\right)-\sin\left(2\pi y_p'(T_{+})+2\pi\omega T_{+}\right))}{\tau_p}dt\\
%&\geq& \int_{T+}^{T} \frac{2a'\cos\left(2\pi y_p'(t)+2\pi\omega t+2\pi y_p'(T_{+})+2\pi\omega T_{+}\right)\sin\left(2\pi y_p'(t)+2\pi\omega t-2\pi y_p'(T_{+})-2\pi\omega T_{+}\right))}{\tau_p}dt\\
%\end{eqnarray*}
%If $W_{+}>0$ then $y_{p}'$ is locally increasing, and the sin in the integral is locally positive. It forces the cosine to be locally nonpositive..........

If $V_{+}>0$ we have $\varepsilon \in [0,1/2]$. And since $V_{i}'$ is decreasing, and $Y_{i}'$ increasing unbounded, there is a moment where $2\pi Y_{i}'(t) + 2\pi\omega t = 2k\pi + \pi = 2\pi Y_{i}'(T_{+})+2\pi\omega T_{+} +\pi/2 - \varepsilon \pi $. At this moment, $V_{i}'$ becomes negative and $Y_{i}'$ becomes decreasing. Since $V_{i}'$ is bounded, it will reach a minimum (since it cannot converges). Denote this time $T_{-}$ and we are in the symmetric case than previously.

We have proved that there exists two sequences $(T_{+}^{n})_{n\in\mathbb{N}}$ and $(T_{-}^{n})_{n\in\mathbb{N}}$ such that $T_{+}^{n} < T_{-}^{n} < T_{+}^{n+1}$ for all $n\in \mathbb{N}$. We can bounded the time $(T_{-}^{n}-T_{+}^{n})$ above and below independently of $n\in \mathbb{N}$ roughly proving that the solution is close to a periodic one. Finally the solution $X_{i}'$ is close to a increasing signal having periodic oscillation around its drift, which is incompatible with an expected diffusive behavior.
\end{proof}

In this particular case of only one sine, we have performed a transformation which leads to an autonomous system, and hard conclusion with only a discrete set of final positions. With more exciting sines the behavior could be different. But -as it is represented in Fig.~\ref{dispersion-dim}- even with more exciting sines we do not obtain in 1D a dispersive behavior as expected. It makes a 1D model very dubious.

But, dispersion of particles is greatly influenced by the dimensionality of the underlying space chosen. Although the dynamic in the one dimensional case is very different from the physic we aim at modeling, we expect that when dimensionality is increased, this behavior will change and be most likely similar to diffusion (see Fig. \ref{dispersion-dim}), as envisioned by the physic, and as described by the models currently in use in the literature. Let us check this assumption in the following section.

\subsection{Higher dimensionality}

It is possible to observe numerically that by increasing the dimensionality to more than one physical dimension (Figs. \ref{fig:stat-xpn-var-2D} and \ref{fig:stat-xpn-var-3D}), the second order moment of \(f_t\) has a \textit{better behavior}, i.e. it increases quite monotonously with time, and the particles do not seem to be overly constrained by the underlying fluid flow. The higher the dimensionality, the better the dispersion of the particles. Indeed, one observes in Fig. \ref{fig:stat-xpn-var-2D} that the dispersion of the particles appears to be much less influenced by the characteristics of the underlying fluid flow than in the 1D case (see Fig. \ref{fig:stat-xpn-var-1D}), and that the third dimensionality brings even more smoothness (see Fig. \ref{fig:stat-xpn-var-3D}). The change of behavior between 2D and 3D can also be partly understood by the addition of new topologies for the three-dimensional stationary points as described in \cite{bec2005clustering}.

Given these results, it seems relevant to keep on pursuing the simulation effort focusing on the three dimensional configuration.

\FloatBarrier
\section{Towards two-way coupled systems}
\label{Towards two-way coupled systems}
The next step towards the modeling of particulate flows is to account for the impact of the disperse phase on the turbulent carrier phase, which has strong implications. Let us consider the empirical measure $\mu_t^N(\bZ) = \frac{1}{N} \sum_{i=1}^N \delta_{\bX_i(t)}\delta_{\bC_i(t)}$ and the following evolution equation
\beq{eq:Liouville_reduced}
 \left\{
   \begin{array}{ccl}
     \ds d\bX_i &=& \bC_i(t) dt,\\
     \\
     \ds d\bC_i &=& \dfrac{\mathbf{u}_f(t,\bX_i)-\bC_i}{\tau_p}dt%=\G\left(t,\bX_1,\bC_1,\dots,\bX_N,\bC_N \right) dt 
   \end{array}
 \right.
 \quad i=1,\dots,N.
\eeq
In a one-way coupled context the gas phase velocity at the particle location does only depend on the particle position itself and is independent of the others particles as they share the same gas phase. In this context, we satisfy the conditions of Theorem \ref{thrm-julie}, i.e. $\G\left(t,\bX_i,\bC_i\right) = \frac{\mathbf{u}_f(t,\bX_i)-C_i}{\tau_p} $. We can thus state a theorem of convergence towards the law of the process.

 In a two-way coupled system, all particles affect the gas phase evolution such that the gas velocity is conditional to the full particle configuration. It can be parametrized by the empirical measure at time $t=0$:
\beq{eq:gas_withparticle}
\mathbf{u}_f(t,x)=\mathbf{u}_f^N(t,x,\mu_0^N[\bZ])
\eeq
In this case, the drag term now depends on all particle history, i.e. $\G^N\left(t,\bX_i,\bC_i ,\mu_0^N[\bZ] \right) = \frac{\mathbf{u}_f^N(t,\bX_i, \mu_t^N[\bZ])-C_i}{\tau_p}$. Now, we are not in the context of Theorem \ref{thrm-julie} anymore. The open question is then to determine if it is possible to characterize a convergence of the particulate system towards a one-particle law:
\begin{align}\label{eq_ft_tw}
\frac{\partial}{\partial t} f_t + \boldsymbol{v} \cdot \nabla _{\boldsymbol{x}} f_t + \nabla _{\boldsymbol{v}}\cdot \left( \G^{lim}f_t \right) = 0 \\
\partial_t\left(\boldsymbol{u}_f\right)+\left(\boldsymbol{u}_f\cdot\nabla\right)\boldsymbol{u}_f - \nu\nabla^2\boldsymbol{u}_f  =  -\frac{1}{\rho_f} \nabla p +\int m_pG^{lim}\text{d}f_t\boldsymbol{v}
\end{align}
where $G^{lim}$ is the forcing of the gas velocity field for a large number of particles, i.e. when the particulate phase behaves as a continuum, and $m_p$ is the mass of each particle.

\subsection{Example of the Burgers equation}

To investigate if there is an Eulerian continuum limit to the two-way problem, we set up a simplified case that considers the 1D Burgers equation on the gas velocity $u$:
\begin{align}
%\dfrac{\partial u}{\partial t}+\dfrac{\partial}{\partial x}\left(\dfrac{u^2}{2}\right)=\nu \dfrac{\partial^2 u}{\partial x^2}+\dfrac{\boldsymbol{F}_{p->g}(t,\boldsymbol{x})}{\rho_f}
\dfrac{\partial u}{\partial t}+\dfrac{\partial}{\partial x}\left(\dfrac{u^2}{2}\right)=\dfrac{F_{p\rightarrow g}(t,x)}{\rho_f}
\end{align}
Giving a meaning to \(F_{p\rightarrow g}\) is not trivial (see \cite{lagoutiere2008simple, aguillon2015riemann,vignal2001}). Here we will use the numerical cells as a regularization for the particle field. The equation is solved using a 1st order finite volume scheme.

\subsubsection{Solution with homogeneous distribution of particles}
First we study the asymptotic limit in which the particles are perfectly uniformly distributed at time $t=0$ at the same velocity. The gas velocity also starts at a uniform velocity. In this limit the Eulerian continuum limit is valid and the particles can be represented by their eulerian equations. We then state that the forcing term in the kinetic equation is $G^{lim}=\frac{u_g-v}{\tau_p}$.
Coupling gas phase and liquid phase equations, we get:
\begin{align}
\dfrac{\partial m_pn_l}{\partial t}+\dfrac{\partial m_pn_lu_l}{\partial x}&=0 \\
\dfrac{\partial m_pn_lu_l}{\partial t}+\dfrac{\partial m_pn_lu_l^2+P_l}{\partial x}&=n_lm_p\dfrac{u-u_l}{\tau_p} \\
\dfrac{\partial u}{\partial t}+\dfrac{\partial}{\partial x}\left(\dfrac{u^2}{2}\right)&=\dfrac{m_pn_l}{\rho_f}\dfrac{u_l-u}{\tau_p}
%\dfrac{\partial u}{\partial t}+\dfrac{\partial}{\partial x}\left(\dfrac{u^2}{2}\right)&=\nu \dfrac{\partial^2 u}{\partial x^2}+\dfrac{m_pn_l}{\rho_f}\dfrac{u_l-u}{\tau_p}
\end{align}
where $m_p$ is the (constant) mass of one particle, $n_l$ the number of particles per unit volume and $P_l$ the pressure of the dispersed phase. In the following we make the assumption of monokinetic disperse phase, i.e. $P_l=0$. The gas density $\rho_f$ is also assumed to be constant. Starting from an homogeneous conditions, we can easily see that the solution will still be invariant by translation at any time and the problem to be solved reduces in the following ODE:
\begin{align}
\dfrac{\text{d}u_l}{\text{d}t}&=\dfrac{u-u_l}{\tau_p} \\
\dfrac{\text{d}u}{\text{d}t}&=\dfrac{m_pn_l}{\rho_f}\dfrac{u_l-u}{\tau_p} 
\end{align}
which solution is:
\begin{align}
u_l(t)&=-\dfrac{1}{\kappa}\left(u^0-\dfrac{\kappa u_l^0+u^0 }{1+\kappa}\right)e^{-\frac{1+\kappa}{\tau_p}t}+\dfrac{\kappa u_l^0+u^0 }{1+\kappa}  \\
u(t)&=\left(u^0-\dfrac{\kappa u_l^0+u^0 }{1+\kappa}\right)e^{-\frac{1+\kappa}{\tau_p}t}+\dfrac{\kappa u_l^0+u^0 }{1+\kappa}
\end{align}
where $\kappa=\frac{m_pn_l}{\rho_f}$. The equilibrium solution is then:
\begin{align}
u_l(t\rightarrow\infty)=u(t\rightarrow\infty)=\dfrac{\kappa u_l^0+u^0 }{1+\kappa}
\end{align}
As a consequence, if we want to study the impact of inhomogeneity of the particulate phase by changing the number of particles but keeping the same physical problem, we need to modify the particle mass $m_p$ accordingly, to keep $\kappa_m=\int \kappa \text{d}x/L_x=N_pm_p$ constant.

\subsection{Particle-laden case with Lagrangian particles}
\begin{figure}
\centering
 	\includegraphics[width=0.9\textwidth]{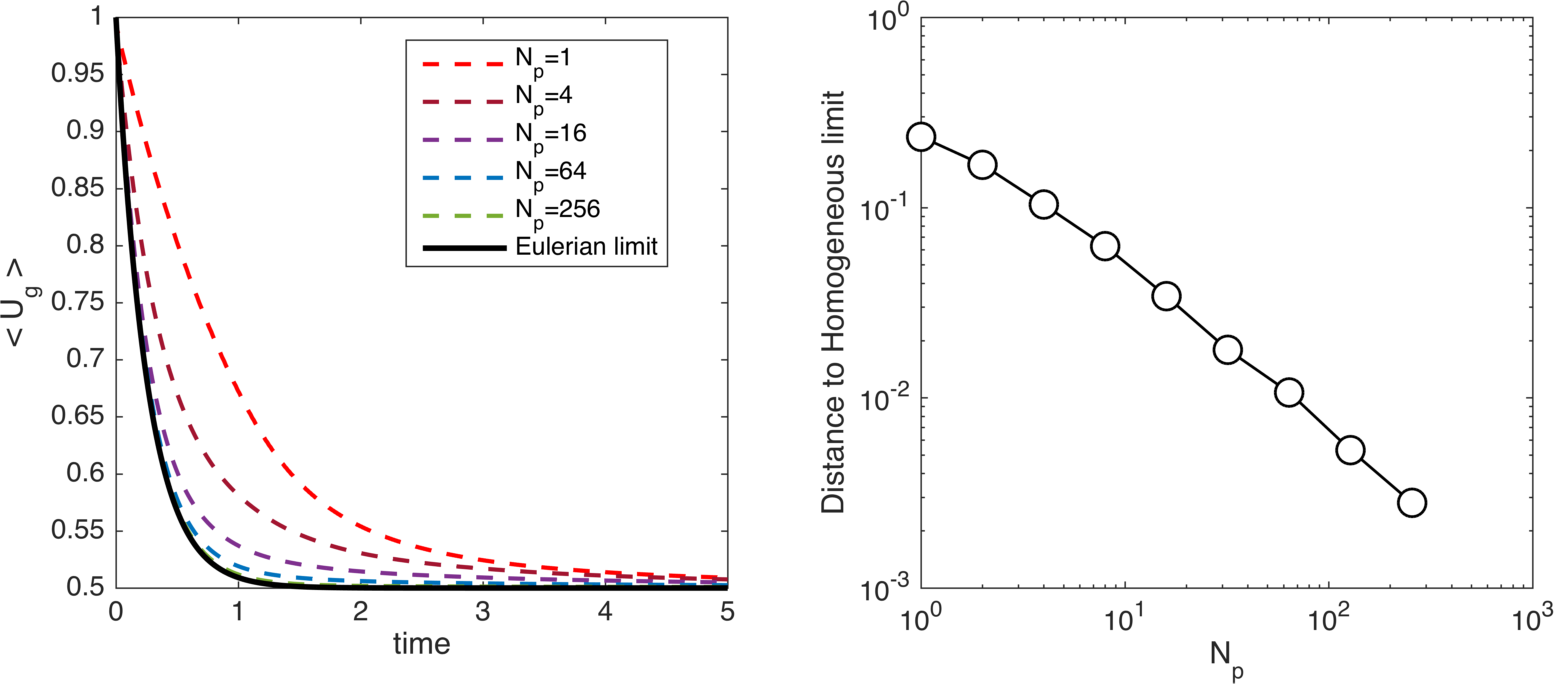}
\caption{1D burgers problem. Left: ensemble-averaged gas velocity for different number of particles from $1$ to $256$, to compare with the homogeneous limit (black). Right: convergence rate of the difference with the homogeneous solution with respect to the number of particles.}
\label{burgers_lag}
\end{figure}

 Knowing the sought continuum limit of the particle system, we now investigate the impact of the number of particles, i.e. the impact of the statistical convergence of the randomly-drawn initial condition. We thus simulate the two-way coupled burgers problem by changing the number of particles from 1 to a large number or particles. In Fig.~\ref{burgers_lag}a, we compare the time evolution of the gas velocity averaged over a large number of realizations of the initial conditions for different numbers of particles at fixed mass loading. We clearly see the convergence of the Lagrangian simulations towards the homogeneous solution, with a convergence rate of order $1$ (see Fig.~\ref{burgers_lag}b). This convergence rate is not affected by the number of cells for numerical discretization and by the addition of physical diffusion in the Burgers equation. So even if we do not have a formal proof in the spirit of Theorem \ref{thrm-julie}, we still have confidence in the existence of a convergence result, and thus of an Eulerian limit description. 

\subsection{Eulerian modeling}
As examplified by the previous test case, the Eulerian representation is still possible for a large number of particles. If we do not have a large number of particles, let say less than one particle per cell, we still have to propose a closure for $G^{lim}$. Moreover, having a statistically-converged NDF $f$ requires to take statistics also on the gas phase velocity. We thus have a two-fold closure problem:
\begin{align}
\frac{\partial}{\partial t} f_t + v \cdot \nabla _x f_t + \nabla _v\cdot \left( G^{lim}f_t \right) = 0, \\
\dfrac{\partial \left<u\right>}{\partial t}+\dfrac{\partial}{\partial x}\left(\dfrac{\left<u\right>^2}{2}\right)=\left<\dfrac{F_{p\rightarrow g}(t,x)}{\rho_f}\right>-\dfrac{\partial}{\partial x}\left(\dfrac{\left<u'\right>^2}{2}\right)
\end{align}
where $<.>$ stands for the ensemble-average over particle realizations which clearly denotes a mean over the initial law of particles.

Here we clearly see that performing an Eulerian simulation sought as an ensemble-average simulation necessarily leads to an ensemble-average on the gas phase. 
Closing the whole system is a tough task outside of the scope of the present work. 
\subsubsection{Closing the equations}
\begin{figure}
\centering
 \subfigure[]          
 	{\includegraphics[height=5cm]{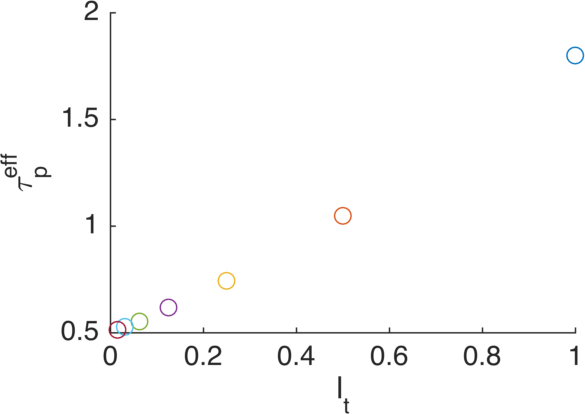} }
 \subfigure[zoom-in]      
 	{\includegraphics[height=5cm]{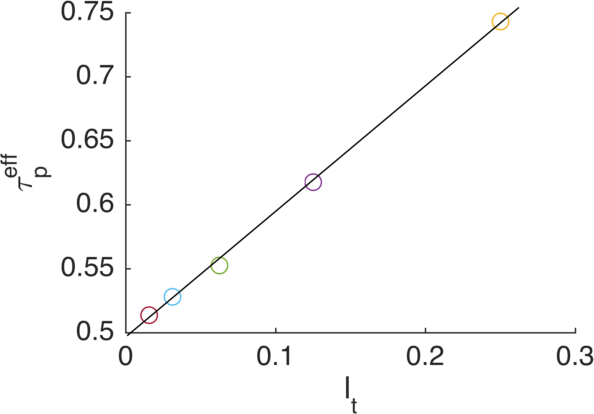}}
\caption{1D burgers problem. Evaluation of an effective relaxation time scale for the Eulerian modeling as a function of the droplet interspace.}
\label{burgers_eul_fit}
\end{figure}

Here we give some insight of possible closures. As results in Fig.~\ref{burgers_lag} clearly shows similar trends but with a different time scale, we propose to investigate the possibility to close the problem using an adapted relaxation time scale $\tau_p^{eff}$:
\begin{align}
\dfrac{\text{d}u_l}{\text{d}t}&=\dfrac{u-u_l}{\tau_p^{eff}} \\
\dfrac{\text{d}u}{\text{d}t}&=\kappa\dfrac{u_l-u}{\tau_p^{eff}} 
\end{align}
In Fig.~\ref{burgers_eul_fit}, we look at the impact of the particle interspace $l_t=1/N_p$ on this effective time scale. We exhibit a linear trend for small $l_t$, which would be helpful to devise closures in a two-way coupled system. The closure for this effective time scale can then be sought as:
$$ \tau_p^{eff}=\tau_p+\alpha l_t$$

The previous example was just to show the impact of the droplet interspace on the solution, and a possible modeling strategy to account for some of the effects. We only focused on the source term, but additional fluxes can also be investigated as possible closures.

\subsubsection{Interpreting existing Eulerian simulations}
Even if we clearly show here that the ensemble average on the particle phase leads to an ensemble-average on the gas phase, i.e. RANS-like statistics, many simulations can be found in the literature in a LES context, which obviously exhibits large scale unsteady behavior. Thus, the question is: what is solved in such simulations? 
A possible and simple interpretation is not to consider this simulation as statistics but as a unique realization of the disperse phase represented in a Eulerian manner. This way, a unique realization of a gas phase will be considered. This turns out to be an Eulerian representation of the empirical measure, which is valid in the sense of the distributions:
\begin{equation}\label{eq_mutn_tw}
\frac{\partial}{\partial t} \mu_t^N + v \cdot \nabla _x \mu_t^N + \nabla _v\cdot \left(\dfrac{u_g-v}{\tau_p}\mu_t^N \right) = 0.
\end{equation}
Taking the moments of this equation and the gas equation, we get:
\begin{align}
\dfrac{\partial m_pn_l^N}{\partial t}+\dfrac{\partial m_pn_l^Nu_l^N}{\partial x}&=0 \\
\dfrac{\partial m_pn_l^Nu_l^N}{\partial t}+\dfrac{\partial m_pn_l^N{u_l^N}^2+P_l^N}{\partial x}&=n_l^Nm_p\dfrac{u-u_l^N}{\tau_p} \\
\dfrac{\partial u}{\partial t}+\dfrac{\partial}{\partial x}\left(\dfrac{u^2}{2}\right)&=\dfrac{m_pn_l^N}{\rho_f}\dfrac{u_l^N-u}{\tau_p}
%\dfrac{\partial u}{\partial t}+\dfrac{\partial}{\partial x}\left(\dfrac{u^2}{2}\right)&=\nu \dfrac{\partial^2 u}{\partial x^2}+\dfrac{m_pn_l}{\rho_f}\dfrac{u_l-u}{\tau_p}
\end{align}
where $n_l^N$ and $u_l^N$ are zeroth and first order moments of the empirical measure, and $P_l^N$ is its pressure. This system of equations is similar to the Eulerian continuum limit, but the difference lies in the initial and boundary conditions: while for the continuum limit, these inputs must be related to the law, here they must randomly drawn as in the case of the Lagrangian particles.

\begin{figure}
\centering
\includegraphics[height=6cm]{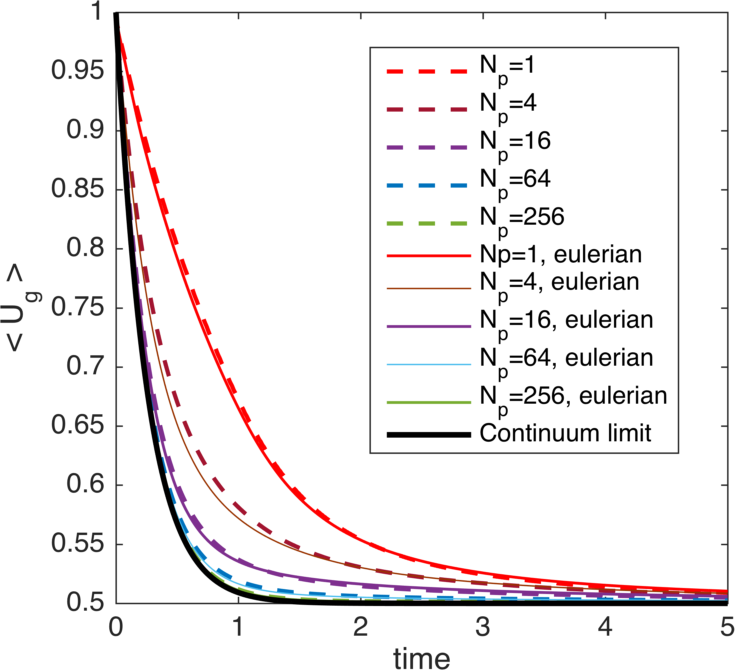} 

\caption{1D burgers problem. Ensemble-averaged gas velocity for different number of particles from $1$ to $256$, with the Lagrangian tracking (dashed lines) and the "empirical" Eulerian moment method (full lines to compare with the homogeneous limit (black).}
\label{burgers_eul_empirical}
\end{figure}

In the case of the 1D burgers periodic problem, solving this system will not take advantage of the spatial invariance of the problem, and we thus have to solve the PDEs. In the following, we will consider a pressureless dynamics, i.e. $P_l^N=0$, and we will use a second order scheme considering the high number density gradients to be resolved.
In Fig.~\ref{burgers_eul_empirical}, we show the results of the gas phase statistics when using this "empirical" Eulerian moment method, demonstrating the ability of such representation to capture the right behavior.

At this point, it is worth to mention that existing LES two-way coupled simulations do not consider a random sampling of the initial/boundary conditions. Instead, they use statistically-converged inputs, leading to an incoherent modeling. It is possible to consider them as regularized simulations in the sense that the initial/boundary conditions has been smeared out enough to lose any random effect.

\section{Conclusions}
In this paper, an exploratory research activity has been started with the aim of statistical and probabilistic modeling of a cloud of particles coupled with a turbulent fluid. Regarding the complexity of this problem, the wide range of expertise of the authors represents an important asset. Here we have set up a common basis to address the issues arising from the context of this work. By investigating all the passing to the limit, we have clarified the main milestones to reach in order to answer our problematic. We have also defined a proper numerical framework to evaluate the modeling approaches and to investigate the statistical properties of our systems of interest. Finally, we have shown the main limitations in two-way coupled system, proposing some possible solutions to overcome them.

\bigskip
\begin{bf}Acknowledgement~:\end{bf}
The financial support by D.G.A. for PhD thesis of D. Mercier is gratefully acknowledged.

%%-----------------------------
%%      your bibliography
\bibliographystyle{abbrv}
\bibliography{proceeding}
%%-----------------------------
\end{document}